\documentclass[11pt,reqno,tbtags]{amsart} 
\usepackage{upref}

\theoremstyle{plain}
\newtheorem{theorem}{Theorem}

\newtheorem{lemma}[theorem]{Lemma}
\newtheorem{corollary}[theorem]{Corollary}
\newtheorem{proposition}[theorem]{Proposition}

\theoremstyle{definition}
\newtheorem*{ack}{Acknowledgement}

\newtheorem{remark}[theorem]{Remark}

\theoremstyle{remark}

\numberwithin{equation}{section}
\numberwithin{theorem}{section}

\newenvironment{romenumerate}{\begin{enumerate}
 }{\end{enumerate}}

\newenvironment{thmxenumerate}{\begin{enumerate}
\setlength{\leftmargin}{0pt}
\setlength{\itemindent}{0pt}
 }{\end{enumerate}}

\newcommand{\refS}[1]{Section ~\ref{#1}}
\newcommand{\refA}[1]{Appendix ~\ref{#1}}
\newcommand{\refT}[1]{Theorem ~\ref{#1}}
\newcommand{\refC}[1]{Corollary ~\ref{#1}}
\newcommand{\refP}[1]{Proposition ~\ref{#1}}
\newcommand{\refL}[1]{Lemma ~\ref{#1}}
\newcommand{\refR}[1]{Remark ~\ref{#1}}
\newcommand{\refand}[2]{\ref{#1} and~\ref{#2}}

\begingroup
  \divide\count255 by 60
  \multiply\count255 by -60
  \advance\count255 by \time
  \ifnum \count255 < 10 \xdef\klockan{\the\count1.0\the\count255}
  \else\xdef\klockan{\the\count1.\the\count255}\fi
\endgroup

\newcommand\nopf{\qed}   

\renewcommand\phi{\varphi}

\renewcommand{\=}{:=}

\newcommand{\aex}{a.e.\spacefactor=1000}
\newcommand\set[1]{\ensuremath{\{#1\}}}

\newcommand\xpar[1]{(#1)}
\newcommand\bigpar[1]{\bigl(#1\bigr)}
\newcommand\Bigpar[1]{\Bigl(#1\Bigr)}
\newcommand\biggpar[1]{\biggl(#1\biggr)}
\def\rompar(#1){\textup(#1\textup)}    
\newcommand\xfrac[2]{#1/#2}
\newcommand\parfrac[2]{\Bigl(\frac{#1}{#2}\Bigr)}

\newcommand\ceil[1]{\lceil#1\rceil}
\newcommand\floor[1]{\lfloor#1\rfloor}

\newcommand\ntoo{\ensuremath{n\to\infty}}
\newcommand\ktoo{\ensuremath{k\to\infty}}

\newcommand\xtoo{\ensuremath{x\to\infty}}

\newcommand\ltoox[1]{\ensuremath{\gl\to#1\infty}}

\newcommand\E{\operatorname{\mathbb E{}}}
\renewcommand\P{\operatorname{\mathbb P{}}}
\newcommand\Var{\operatorname{Var}}

\newcommand\Po{\operatorname{Po}}
\newcommand\Bi{\operatorname{Bi}}
\newcommand\ett{\boldsymbol1}

\newcommand\dd{\,d}

\newcommand\intoo{\int_0^\infty}
\newcommand\intoe{\int_0^\eps}
\newcommand\inteo{\int_\eps^\infty}
\newcommand\limn{\lim_{n\to\infty}}

\newcommand\qww{^{2/3}}
\newcommand\qw{^{1/3}}
\newcommand\qq{^{1/2}}
\newcommand\qqi{^{-1/2}}
\newcommand\qi{^{-1}}
\newcommand\qii{^{-2}}
\newcommand\qiw{^{-1/3}}
\newcommand\qiww{^{-2/3}}

\newcommand\pfitem[1]{\par(#1):}

\newcommand\ga{\alpha}

\newcommand\gd{\delta}

\newcommand\gl{\lambda}
\newcommand{\lam}{\lambda}
\newcommand\gL{\Lambda}
\newcommand\eps{\varepsilon}
\newcommand\cL{\mathcal L}
\newcommand\cA{B}
\newcommand\fN{\mathfrak N}
\newcommand\fS{\mathfrak S}

\newcommand\iid{i.i.d.\spacefactor=1000}     
\newcommand\ie{i.e.\spacefactor=1000}
\newcommand\eg{e.g.\spacefactor=1000}

\newcommand\cf{cf.\spacefactor=1000}
\newcommand{\as}{a.s.\spacefactor=1000}

\newcommand\dto{\overset{\mathrm{d}}{\to}}
\newcommand\pto{\overset{\mathrm{p}}{\to}}
\newcommand\asto{\overset{\mathrm{a.s.}}{\to}}
\newcommand\eqd{\overset{\mathrm{d}}{=}}
\newcommand\psim{\sim_p}

\newcommand\bbR{\mathbb R}
\newcommand\bbN{\mathbb N}
\newcommand\bbNx{\mathbb N^*}

\newcommand\cC{\mathcal C}

\newcommand\Cc{C_c}
\newcommand\xini{\xi_{ni}}

\newcommand\fall[2]{#1^{\underline{#2}}}

\newcommand{\G}{G}
\newcommand\gnx[1]{\ensuremath{\G(n,#1)}}
\newcommand\gnp{\gnx{p}} 
\newcommand\gnm{\gnx{m}} 

\newcommand\ze{Z_\eps}
\newcommand\zne{Z_{n\eps}}
\newcommand\znex{Z_{n\eps}'}
\newcommand\we{W_\eps}
\newcommand\wne{W_{n\eps}}
\newcommand\wl{w_\ell}
\newcommand\xil{\Xi^\ell}
\newcommand\xilx[1]{\Xi^{#1}}
\newcommand\Xila{\Xi\gll}
\newcommand\Xilas{\Xi\glls}
\newcommand\xila{\xi\gll}
\newcommand\xilai{\xila_i}
\newcommand\Xix{\Xi^*}
\newcommand\gll{^{(\lam)}}
\newcommand\glls{^{(\lam-s)}}
\newcommand\gllx{^{(\lam-x)}}
\newcommand\gxx[1]{^{(#1)}}
\newcommand\LL{L}
\newcommand\ue{U_{\eps}}
\newcommand\uel{\ue(\gl)}
\newcommand\glx{\gamma}

\newcommand\sumke{\sum_{k=1}^{\eps n^{2/3}}}

\newcommand\sumje{\sum_{j=1}^{\eps n^{2/3}}}
\newcommand\sumkab{\sum_{k=a n\qww}^{b n^{2/3}}}
\newcommand\sumjab{\sum_{j=a n\qww}^{b n^{2/3}}}
\newcommand\iab{[a n\qww,b n^{2/3}]}
\newcommand\iie{[1,\eps n^{2/3}]}

\newcommand\ekk{E_{k,k}}
\newcommand\ekj{E_{k,j}}
\newcommand\ek{E_{k}}
\newcommand\ej{E_{j}}
\newcommand\tkx{t_k^*}
\newcommand\gex{\ge\!}
\newcommand\gexx{\ge\!\!}
\newcommand{\nc}{n_{c}}
\newcommand\pkp{P_{k,p}}
\newcommand\pkk[1]{P(n; k_1,\dots,k_j;#1)}
\newcommand\hla[1]{h\gll_{#1}}
\newcommand\al{a_\gl}
\newcommand\alx{a_{\gl'}}
\newcommand\dtv{d_{\mathrm{TV}}}
\newcommand\phil{\phi_\gl}

\newcommand\lhs{left hand side}
\newcommand\rhs{right hand side}

\newcommand\CS{Cauchy--Schwarz}

\begin{document}

\title[{A point process 
of the random graph evolution}]
{A point process describing 
the component sizes 
in the critical window
of the random graph evolution}
\author{Svante Janson}
\address{Department of Mathematics, Upp\-sala University, PO Box 480,
S-751 06 Upp\-sala, Sweden}
\email{svante.janson@math.uu.se}
\urladdr{http://www.math.uu.se/\~{}svante/}

\author{Joel Spencer}
\address{Joel Spencer,
Courant Institute,
251 Mercer St.,
New York, NY 10012, USA}
\email{spencer@cs.nyu.edu}
\urladdr{http://www.cs.nyu.edu/cs/faculty/spencer/}


\date{May 24, 2005} 

\begin{abstract}
We study
a point process describing the asymptotic behavior of 
sizes of the largest components of the random graph \gnp{}
in the critical window $p=n^{-1}+\gl n^{-4/3}$.
In particular, we show that this point process has a surprising rigidity.
Fluctuations in the large values will be balanced by opposite
fluctuations in the small values such that the sum of the values
larger than a small $\eps$ is almost constant.
\end{abstract}


\maketitle

\section{Introduction}
We consider the asymptotic behavior of the component sizes in the
random graph \gnp, where throughout this paper $p=n^{-1}+\gl n^{-4/3}$
for some fixed $\gl$ with $-\infty<\gl<\infty$. 
It is well-known 
that this is the critical window of $p$
where the ``phase transition'' occurs.
It is further well-known that, for the $p$ we consider, the largest components
are of order $n^{2/3}$. We therefore scale by this factor;
if the components are $\cC_1,\cC_2,\dots,\cC_r$, in order of
decreasing size, say,
and $|\cC_i|$ is the size (order) of $\cC_i$, 
we define $\xi_{ni}$ to be $n^{-2/3}|\cC_i|$
and consider the random set
$\Xi_n:=\set{\xi_{ni}}_{i=1}^r$ as a point process on $(0,\infty)$
or $(0,\infty]$. See
\refA{S:point} for some technical background and note that it is
convenient to define the point process formally as a random measure
with point masses at the points $\xi_{ni}$; we will sometimes use this
formalism, writing for example $\Xi_n[a,b]$ for the number of points
in $[a,b]$, but we will also speak
(and think) of
point processes as random sets.

It follows immediately from Aldous \cite[Corollary 2]{Aldous}, see
\refL{L:point2}, that as \ntoo, 
the point processes $\Xi_n$ converge in distribution to some point
process $\Xila=\set{\xilai}$ on $(0,\infty]$ (in the vague topology on
$(0,\infty]$, see \refA{S:point}); this also follows from 
a minor extension
of results in {\L}uczak, Pittel and Wierman \cite{LPW},
see Janson, {\L}uczak and Ruci\'nski \cite[Theorem 5.20]{JLR}.
Aldous \cite{Aldous} further gave a description of the limiting process $\Xila$
as the set of lengths of excursions of a certain reflected Brownian
motion with parabolic drift, defined as $B^\gl(s):=W^\gl(s)-\min_{0\le u\le s}
W^\gl(u)$, $s\ge0$, where $W^\gl(s)=W(s)+\gl s-s^2/2$ for a standard Brownian
motion $W$. 

We will usually keep $\gl$ fixed, and will then often omit it from the
notation, thus writing $\Xi=\Xila$ and $\xi_i=\xilai$.
Conversely, we may write $\Xi_{n,p}$ when necessary.
Note that we may regard $(\Xila)_\gl$ as a
stochastic process indexed by $\gl\in(-\infty,\infty)$; this is the
standard multiplicative coalescent as constructed by Aldous
\cite{Aldous}, except that the variables $\Xila$ are represented 
as point processes while Aldous uses the equivalent
representation as sequences $(\xi_i)_1^\infty$;
\cf{} \refL{L:point2}, although Aldous uses a stronger topology.

The aim of this paper is to study the limiting point process $\Xi$.
The number of components in \gnp{} tends to infinity (in probability)
as \ntoo, so we expect an infinite number of points $\xi_i$ in $\Xi$.
Moreover, 
if we say that the \emph{weight} of a point $x$ is $x$,
the total weight of $\Xi_n$ is $\sum_i\xi_{ni}=n^{1/3}$, so we expect the
total weight of $\Xi$, \ie{} $\sum \xi_i=\int x\,d\Xi$, to be
infinite a.s.; indeed, this  
is a simple consequence of \refT{T1}.
(Still, we caution that results on the limiting process $\Xi$ do not
automatically follow from results on the discrete $G(n,p)$.  Had we,
for example, chosen the ``wrong'' parameterization 
$\xini=n^{-0.7}|\cC_i|$ then $\Xi$ would be almost surely empty.)

Our main result is the following.
We also give in later sections
various other results; 
several of them have been more or less well-known for a long time, but
perhaps not published previously in this form. 

\begin{theorem}\label{T1}
Let $-\infty<\gl<\infty$, and let $\Xi$ be the limiting point process
defined above.
Let $\ze\=\sum_{\xi_i\ge\eps}\xi_i=\int_\eps^\infty x\,d\Xi(x)$ be the total
weight of all points in $\Xi$ that are at least $\eps$. Then, as $\eps\to0$,
\begin{align}\label{t1a}
\E \ze
&=\parfrac{2}{\pi}^{1/2}\eps^{-1/2}+\gl+
(2\pi)^{-1/2}\lam^2\eps^{1/2}+ O(\eps)
\\
\intertext{and}
\label{t1b}
\Var\ze
&=(2/\pi)\qq\eps^{1/2}+O(\eps).
\end{align}
In particular, $\E\ze\to\infty$ and
$\Var\ze\to0$ as $\eps\to0$.
\end{theorem}

We will also give an exact, but more complicated, formulas for $\E\ze$
in \refC{C1}
and $\Var\ze$ in \refC{Cvar}.
It seems non-trivial to obtain the asymptotics above from these
formulas. 

Thus, as $\eps\to0$, the variables $\ze$ tend to infinity, but 
they become more and more concentrated about their mean; hence,
the random fluctuations disappear in the limit.
In other words, the process $\Xi$ is very rigid, and any random
fluctuation in the weights of the largest points has
to be exactly balanced by opposite fluctuations in the weights of
smaller points; this will be seen again in \refS{Spalm} where we
consider the Palm distributions.
Note that, because of the scaling, this is a non-trivial
result in contrast to 
the corresponding fact that $\Xi_n$ has a constant total weight
$n^{1/3}$.  
Note also that this is very far from the behaviour of a Poisson process.

We will prove \refT{T1} by two different methods, both classical,
each giving a 
partial result only, in Sections \refand{Spf1}{Spf2}.

In contrast to \refT{T1}, the number of points $\gex\eps$ in $\Xi$,
\ie{} $\Xi[\eps,\infty)$, is not sharply concentrated. 

\begin{theorem}\label{T2}
Let $\we\=\Xi[\eps,\infty)$ be the number
of points in $\Xi$ that are at least $\eps$. 
Then, as $\eps\to0$,
\begin{align}\label{t2a}
\E \we
&=
\parfrac{2}{9\pi}^{1/2}\eps^{-3/2}
-(2\pi)^{-1/2}\lam^2\eps^{-1/2}+ \frac14\ln(1/\eps)+O(1)
\\
\intertext{and}
\label{t2b}
\Var\we
&=\parfrac{2}{9\pi}^{1/2}\eps^{-3/2}+O(\eps\qi) 
\sim\E\we.
\end{align}
\end{theorem}

\begin{remark}
  It seems likely that $\we$ is almost Poisson distributed, in the
  sense that its total variation distance to a Poisson distribution
  with the same mean tends to 0 as  $\eps\to0$, but we leave this
  as an open problem. If this holds, it would immediately
  imply asymptotic normality of $\we$.
\end{remark}

The main interest in \refT{T1} comes from the fact that $\ze$
approximatively describes the large component sizes in \gnp{} for
large $n$. We formalize this in 
the following intuitively obvious result; see \refS{Smoments} for a formal
verification of the technicalities.

\begin{proposition}\label{P1}
Let 
$\zne=\sum_{\xini\ge\eps}\xini=\int_\eps^\infty x\,d\Xi_n(x)$
be the total
weight of all points in $\Xi_n$ that are at least $\eps$;
thus $\zne$ equals $n\qiww$ times the total size of all components
$\gex\eps n\qww$ in \gnp.
For every fixed $\eps>0$, as \ntoo, $\zne\dto\ze$
 with convergence of all moments, \ie,
for every $q\ge0$, $\E\zne^q\to\E\ze^q$.

The same holds for $\znex\=\sum_{\xini>\eps}\xini$.

Similarly, if\/ $\wne\=\#\set{i:\xini\ge\eps}=\Xi_n[\eps,\infty)$, 
the number of components $\gex\eps n\qww$ in \gnp, then 
$\wne\dto\we$, with convergence of all moments.
\end{proposition}

We introduce some more notation.
Let
$X_n(k)$ denote the number of components with
$k$ vertices in the random graph \gnp, and let
$Y_n(k):=kX_n(k)$, the total number of vertices in these components.
We further define 
$X_n(I):=\sum_{k\in I}X_n(k)$ 
and
$Y_n(I):=\sum_{k\in I}Y_n(k)$ 
for an interval $I$.
(For simplicity, we omit $p$ from the notation.)
Thus 
$\zne=n\qiww Y_n[\eps n\qww,\infty)$,
$\znex=n\qiww Y_n(\eps n\qww,\infty)$,
and
$\wne=X_n[\eps n\qww,\infty)$.
We denote falling factorials by $\fall nk\=n\dotsm(n-k+1)$.

\begin{remark}
Although we keep $\gl$ fixed for simplicity, it is easy to see
(\eg{} using the monotonicity of \gnp{} in $p$)
that the \refP{P1}  holds also for a sequence $\gl_n\to\gl$. 
Moreover, as a consequence of this and monotonicity,
the same holds for the random graph $G(n,m)$ with a deterministic number 
$m=n/2+(\gl+o(1))n\qww/2$ 
edges.
\end{remark}

\begin{ack}
  We thank David Aldous for interesting discussions.
\end{ack}

\section{First partial proof of \refT{T1}}\label{Spf1}

If $\mu$ is a probability distribution on the non-negative integers,
let $T(\mu)$ denote the (random) total progeny of a
Galton--Watson process with offspring distribution $\mu$, starting
with one initial particle.
(Thus $T(\mu)\in\set{1,2,\dots,\infty}$.)

\begin{lemma}\label{L1}
Let $\gl\in(-\infty,\infty)$ and $\eps>0$ be fixed.
\begin{thmxenumerate}
\item
The limit
$\uel\=\lim_{\ntoo}
n\qw \P\bigpar{T(\Po(1+\gl n\qiw))\ge\eps n\qww}$
exists, and
\begin{equation}\label{ue1}
  \uel=2\max(\gl,0)+
\int_{\eps}^\infty(2\pi)\qqi x^{-3/2}e^{-\gl^2x/2} \dd x.
\end{equation}
\item
More generally, if $\gl_n\to\gl$, then 
$$n\qw \P\bigpar{T(\Po(1+\gl_n n\qiw))\ge\eps n\qww}\to\uel.$$
\item
Moreover, for any fixed $\gd\in\bbR$ and any sequence $\gd_n\to\gd$,
\begin{equation*}
n\qw \P\bigpar{T(\Bi(\floor{n-\gd_n n\qww},n\qi+\gl n^{-4/3}))\ge\eps n\qww}
\to\ue(\gl-\gd).
\end{equation*}
\end{thmxenumerate}
\end{lemma}

\begin{proof}
\pfitem{i}
First consider the possibility of an infinite total progeny. 
By elementary branching process theory,
if $q\=\P(T(\Po(\glx)=\infty))$, then $q=0$ for $\glx\le1$, while $q>0$
for $\glx>1$, and then $1-q=e^{-q\glx}$, or
$\glx=-\ln(1-q)/q=1+q/2+O(q^2)$.
It follows that if $\glx\to1$ with $\glx>1$, then $q\to0$ and 
$q\sim 2(\glx-1)$.
Consequently, as \ntoo, for any real $\gl$,
\begin{equation}\label{emmsan}
n^{1/3}  \P\bigpar{T(\Po(1+\gl n\qiw))=\infty}\to 2\max(\gl,0).
\end{equation}

Next, consider a finite total progeny $\ge\eps n\qww$.
By Otter \cite{Otter}, see also Pitman \cite{Pitman},
  \begin{equation}\label{otter}
\P\bigpar{T(\mu)=k}	
=\frac1k \P(S_k=k-1),
\qquad 1\le k<\infty,
  \end{equation}
where $S_k$ is the sum of $k$  independent random variables with the
distribution $\mu$.
In particular, for a Poisson distribution, using Stirling's formula,
  \begin{equation}\label{emma}
	\begin{split}
\P\bigpar{T(\Po&(1+\gl n\qiw))=k}	
=\frac1k \P\bigpar{\Po(k(1+\gl n\qiw))=k-1}
\\&
=\frac{k^{k-1}(1+\gl n\qiw)^{k-1}}{k!} e^{-k(1+\gl n\qiw)}
\\&
=(2\pi)\qqi k^{-3/2}(1+\gl n\qiw)^{k}e^{-k\gl n\qiw}
\bigpar{1+O(|\gl|n\qiw+k\qi)}.	  
	\end{split}
  \end{equation}
We have
\begin{equation*}
 \ln(1+\gl n\qiw)
= \gl n\qiw-\tfrac12(\gl n\qiw)^2+O(|\gl n\qiw|^3).
\end{equation*}
Hence, for $\eps n\qww\le k< a_n\=(n^{1/3}|\gl|\qi)^{5/2}$
(with $a_n=\infty$ when $\gl=0$),
\eqref{emma} yields
  \begin{equation}\label{emma1}
	\begin{split}
\P\bigpar{T(\Po&(1+\gl n\qiw))=k}	
=(2\pi)\qqi k^{-3/2}e^{-\frac12\gl^2 kn\qiww}
\bigpar{1+O(n^{-1/6})}.	  
	\end{split}
  \end{equation}
Moreover, 
assuming that $n$ is so large that $|\gl n\qiw|<1/2$,
\begin{equation*}
  \ln(1+\gl n\qiw)
\le \gl n\qiw-\tfrac12(\gl n\qiw)^2+\tfrac23|\gl n\qiw|^3
\le \gl n\qiw-\tfrac16\gl^2 n\qiww
\end{equation*}
and thus, by \eqref{emma},
\begin{equation}\label{emma2}
\P\bigpar{T(\Po(1+\gl n\qiw))=k}	
=
O\Bigpar{k^{-3/2}e^{-\frac16\gl^2 k n\qiww}}.
\end{equation}
Summing over $k\in[\eps n\qww,a_n)$ we find, by \eqref{emma1},
\eqref{emma2} and dominated convergence, 
  \begin{equation}\label{samuel}
	\begin{split}
\hskip6em&\hskip-6em
n^{1/3}\P\bigpar{T(\Po(1+\gl n\qiw))\in[\eps n\qww,a_n)}	
\\&
=n^{1/3} \int_{\ceil{\eps n\qww}}^{\ceil{a_n}}
 \P\bigpar{T(\Po(1+\gl n\qiw))=\floor{x}}\dd x
\\&
= \int_{n\qiww\ceil{\eps n\qww}}^{n\qiww \ceil{a_n}}
n \P\bigpar{T(\Po(1+\gl n\qiw))=\floor{n\qww x}}\dd x
\\&
\to
\int_{\eps}^\infty
(2\pi)\qqi x^{-3/2}e^{-\gl^2x/2} \dd x.
	\end{split}
  \end{equation}
Furthermore, the sum over $k\ge a_n$ is exponentially small by \eqref{emma2}.
Hence, the result follows by \eqref{emmsan} and \eqref{samuel}.
\pfitem{ii}
By the same proof as (i), or by (i) and monotonicity.
\pfitem{iii}
One could use \eqref{otter} and argue as above, but we will instead
use a Poisson approximation.
For any $N$ and $p$, we have the bound on the total variation distance
\begin{equation*}
  d_{TV}\bigpar{\Bi(N,p),\Po(Np)} \le p,
\end{equation*}
see \eg{} \cite[Theorem 2.M]{SJI}.
Hence, using a maximal coupling of $\Bi(N,p)$ and $\Po(Np)$ in each
family, we can couple the Galton--Watson processes with offspring
distributions $\Bi(N,p)$ and $\Po(Np)$ such that the probability that
they differ before they have reached at least $n^{1/3}$ individuals is
at most $n^{1/3}p$; furthermore, conditioned on both reaching
$n^{1/3}$ together, and being equal so far, the probability that 
they differ before they have reached at least $\eps n\qww$ individuals 
is at most $\eps n\qww p$. Hence,
\begin{multline}
  \label{david}
\bigl|
\P\bigpar{T(\Bi(N,p)\ge\eps n\qww)} -\P\bigpar{T(\Po(Np)\ge\eps n\qww)}
\bigr|
\\
\le 
n^{1/3}p+\P\bigpar{T(\Po(Np)\ge n^{1/3}} \eps n\qww p.
\end{multline}
Now, let $N=\floor{n-\gd n\qww}$ and $p=n\qi+\gl n^{-4/3}$ and let
\ntoo.
Then $Np\to1$, and thus, for each fixed $M$ and $n>M^3$, 
\begin{equation*}
  \P\bigpar{T(\Po(Np)\ge n^{1/3}}
\le 
  \P\bigpar{T(\Po(Np)\ge M}
\to
  \P\bigpar{T(\Po(1)\ge M}.
\end{equation*}
Since $T(\Po(1))$ is finite a.s., the latter probability tends to 0 as
$M\to\infty$, and it follows that
$  \P\bigpar{T(\Po(Np)\ge n^{1/3}}\to0$. 
Consequently, the \rhs{} of \eqref{david} is 
$O(n\qiww)+o(1)\cdot O(n\qiw)=o(n\qiw)$. The result follows from
\eqref{david} and (ii), since
\begin{equation*}
  Np=(n-\gd n\qww+O(1))(n\qi+\gl n^{-4/3})
=1+\gl_n n\qiw,
\end{equation*}
with $\gl_n=\gl-\gd+O(n\qiw)\to\gl-\gd$.
\end{proof}

We give alternative formulas for $\uel$ defined in \refL{L1}.

\begin{lemma}\label{L2}
Let $-\infty<\gl<\infty$ and $\eps>0$. Then  
\begin{align}
\uel
&=
\parfrac{2}{\pi}^{1/2}\eps^{-1/2}+\gl+
\int_0^{\eps}(2\pi)\qqi x^{-3/2}\bigpar{1-e^{-\gl^2x/2}} \dd x
\label{ue2}
\\&
=
\parfrac{2}{\pi}^{1/2}\eps^{-1/2}+\gl+
(2\pi)^{-1/2}\lam^2\eps^{1/2}+ O(\gl^4\eps^{3/2}).
\label{ue3}
\end{align}
\end{lemma}

\begin{proof}
First note that in the case $\gl=0$, \eqref{ue1} yields
\begin{equation}\label{ue0}
\ue(0)=	
\int_{\eps}^\infty(2\pi)\qqi x^{-3/2}\dd x
=
(2/\pi)\qq \eps^{-1/2}.
\end{equation}
Since $2\max(\gl,0)=\gl+|\gl|$, \eqref{ue1} further yields
\begin{equation}\label{ue4}
\uel-\ue(0)=\gl+|\gl|-
\int_{\eps}^\infty(2\pi)\qqi x^{-3/2}\bigpar{1-e^{-\gl^2x/2}} \dd x.
\end{equation}
Now, for $\gl\neq0$,
by change of variables and a standard integration by parts,
\begin{equation*}
  \begin{split}
\int_{0}^\infty x^{-3/2}\bigpar{1-e^{-\gl^2x/2}} \dd x
&
=(\gl^2/2)\qq
\int_{0}^\infty y^{-3/2}\bigpar{1-e^{-y}} \dd y
\\&
=(\gl^2/2)\qq \,2\,\Gamma(1/2)
=(2\pi)\qq|\gl|,
  \end{split}
\end{equation*}
and thus \eqref{ue4} yields
\begin{equation}\label{ue5}
\uel=\ue(0)+\gl+
\int_0^{\eps}(2\pi)\qqi x^{-3/2}\bigpar{1-e^{-\gl^2x/2}} \dd x.
\end{equation}
This proves \eqref{ue2}, and \eqref{ue3} follows by the expansion
$1-e^{-\gl^2x/2}=\gl^2x/2+O(\gl^4x^2)$.
\end{proof}

\begin{remark}
  Expression \eqref{ue1} might lead the casual reader to
suppose that $\lambda=0$ was somehow special.  The
equivalent expression \eqref{ue2}, however, shows that
$\uel$ is a smooth function of $\lambda$.  This corresponds
to the  generally held belief that there can be no further
refinements of the critical window, that no value of $\lambda$
is more special than any other, 
and that natural functions
vary smoothly with $\lambda$.
\end{remark}

Returning to the random graphs, note that given the graph \gnp,
the probability that a random vertex belongs to a component of size
at least $\eps n\qww$ is $n\qi Y[\eps n\qww,\infty)=n\qiw\zne$.
Taking expectations we see that 
$\E\zne$ equals $n\qw$ times the probability that a given vertex $v$
belongs to a component of size at least $\eps n\qww$ in \gnp.
We explore the component containg the given vertex by the standard
breadth-first search. 
In this search, we explore first the neighbours of $v$, then their
neighbours, and so on, see \eg{} \cite{Spencer} or \cite[Section 5.2]{JLR}. 
When we explore the neighbours of a vertex, we
find $\Bi(n-m,p)$ new vertices in the component, where $m$ is the
number of vertices found so far.
Thus, 
the process is dominated by 
a Galton--Watson process with offspring distribution $\Bi(n,p)$, 
and, if we stop when we reach $\eps n\qww$
vertices, dominates 
a Galton--Watson process with offspring distribution 
$\Bi(\floor{n-\eps n\qww},p)$; hence, the probability that we find at
least $\eps n\qww$ vertices in the component lies between the
probabilities that these Galton--Watson processes have a total progeny
of at least $\eps n\qww$. Consequently,
\begin{multline*}
\P\bigpar{T(\Bi(\floor{n-\eps n\qww},n\qi+\gl n^{-4/3}))\ge\eps n\qww}
\\
\le
n\qiw \E\zne
\le
\P\bigpar{T(\Bi(n,n\qi+\gl n^{-4/3}))\ge\eps n\qww}.
\end{multline*}
By \refL{L1} and \refP{P1}, this yields
\begin{equation*}
\ue(\gl-\eps)
\le
\E\ze
\le
\uel,
\end{equation*}
and \eqref{t1a} follows by \refL{L2}.

It seems more difficult to estimate $\Var\zne$ by this method, and we
will use another approach in \refS{Spf2}.

\section{Complexity}

The \emph{complexity} $c(G)$ of a graph $G$ with $v$ vertices and $e$
edges is defined by $c(G)\=e-v+1$. Thus the complexity is $0$ for
trees, $1$ for unicyclic connected graphs, and $\ge2$ otherwise.
We say that a connected graph with complexity $\ge2$ is \emph{complex}.

We can refine the point processes $\Xi_n$ and $\Xi$ by considering the
complexities of the components. We can think of this as giving each
point in the processes a label; a point in $\Xi_n$ is labelled
by the complexity of the corresponding component.
Formally, we can think of the labelled versions, $\Xix_n$ and $\Xix$,
say, as point processes on the space $(0,\infty]\times \bbN$, or
better $(0,\infty]\times \bbNx$, where $\bbN=\set{0,1,\dots}$ and 
$\bbNx$ is the compact space $\bbN\cup\set{\infty}$.
The results by Aldous \cite[Corollary 2]{Aldous} and 
{\L}uczak, Pittel and Wierman \cite{LPW} referred to above actually
consider the complexity too, and show that $\Xix_n\dto\Xix$ as \ntoo,
for a suitable labelling $\Xix$ of $\Xi$. 
Aldous \cite[Corollary 2]{Aldous} describes $\Xix$ by the 
process $B^\gl$ defined above: introduce a process of marks on
$(0,\infty)$ that, given $B^\gl$, 
is a Poisson process with intensity $B^\gl(s)\dd s$; then,
as said above, the points $\xi_i$ are
the lengths of the excursions of $B^\gl$, and each excursion is
labelled with the number of marks inside it. In other words, 
given $B^\gl$,
each point $\xi_i$ gets a label that has a Poisson distribution whose
mean is the area under the corresponding excursion, and different
points are labelled independently.
We will give another description in \refT{Tcomp} below.

We let, for $\ell\ge0$, $\xil_n$ be the subset \set{\xi_{ni}:c(\cC_i)=\ell}
of $\Xi_n$ of points with labels $\ell$, \ie{} the set of scaled sizes
of components of \gnp{} with complexity $\ell$. Similarly, let $\xil$
be the subset of $\Xi$ of points with labels $\ell$. Since
$\Xix_n\dto\Xix$, we have $\xil_n\dto\xil$ for every $\ell$.

Let $C(k,\ell)$ be the number of connected graphs with 
complexity $\ell$
on $k$ (labelled) vertices (they thus have $k+\ell-1$ edges).
Thus $C(k,0)$ is the number of trees, and by Cayley's theorem,
$C(k,0)=k^{k-2}$.
More generally, Wright \cite{Wright} proved that for every fixed
$\ell$
\begin{equation}\label{wright}
  C(k,\ell)\sim \wl k^{k+3\ell/2-2}
\qquad\text{as \ktoo},
\end{equation}
for some constants $\wl$, for which Wright \cite{Wright} gave a
recursion formula. (See also \cite[\S8]{giant} and the references
there.)
We have $w_0=1$ and $w_1=\sqrt{\pi/8}$.
It was shown in \cite{Spencer} that 
\begin{equation}
  \wl=\frac{\E \LL^\ell}{\ell!},
\qquad \ell\ge0
\end{equation}
where $\LL$ is the area under a normalized Brownian excursion.
If we introduce the moment generating function $\Psi$ of $\LL$, we
thus have
\begin{equation}
  \label{mgf}
\Psi(t)=\E e^{t \LL}=\sum_{\ell=0}^\infty \wl t^\ell.
\end{equation}
The moments $\E\LL^\ell$ and the moment generating function $\Psi$ had
earlier been studied by Louchard \cite{Lou0,Lou}. 
Note that $\Psi(t)$ is finite for all $t>0$ (and thus \eqref{mgf}
holds for all complex $t$); 
indeed, as remarked in \cite[Remark 3.1]{SJISE} (where $\xi=2\LL$),
it follows from the well-known asymptotics for $\wl$, see \eg{} 
\cite[\S8]{giant} and \cite[Theorem 3.3 and (3.8)]{SJWiener},
that
$\E\LL^\ell\sim\sqrt{18}\,\ell\,(12e)^{-\ell/2}\ell^{\ell/2}$
as $\ell\to\infty$,
and thus \cite[Lemma 4.1(ii)]{SJISE} implies, \cf{} \cite[Remark 4.9]{SJISE},
\begin{equation}
\label{jesper}
\Psi(t) 
\sim
\tfrac12{t^2}e^{\xfrac{t^2}{24}}
\qquad \text{as }
t\to+\infty.
\end{equation}

We now can state the result describing $\Xix$.
For $x>0$, let $P_x$ be the distribution on $\bbN$ given by
\begin{equation}\label{pxl}
  P_x(\ell)=\frac{\wl x^{3\ell/2}}{\Psi\bigpar{x^{3/2}}},
\qquad \ell\ge0.
\end{equation}

\begin{theorem}
  \label{Tcomp}
The point process $\Xix$ on $(0,\infty]\times\bbNx$ can be obtained
from $\Xi$ by independently giving each point $\xi_i\in\Xi$ a random
label with the distribution $P_{\xi_i}$.
\end{theorem}

\begin{proof}
Conditioned on the vertex sets of the components of \gnp,
the internal structures of the components are independent.
Moreover, a component of order $k$ is distributed as $G(k,p)$
conditioned on being connected.
Let $\pkp(\ell)$ be the probability that such a component has
complexity $\ell$. The probability that $G(k,p)$ is connected and has
complexity $\ell$ is $C(k,\ell)p^{k+\ell-1}(1-p)^{\binom k2-k-\ell+1}$ 
and thus
\begin{equation}
  \label{pkp}
\pkp(\ell) 
= C(k,\ell)\parfrac{p}{1-p}^\ell
\Bigm/
\sum_{\ell=0}^\infty
C(k,\ell)\parfrac{p}{1-p}^\ell.
\end{equation}
Consequently, the labelled process $\Xix_n$ can be obtained from
$\Xi_n$ by giving the points $\xini$ labels independently, such that
the label of a point $x$ has the distribution $P_{xn\qww,p}$ given by
\eqref{pkp}. 

By Bollob\'as \cite[Theorem V.20]{Bollobas}, there exists a constant
$c>0$ such that, for all $k$ and $\ell$,
\begin{equation}\label{boll}
  C(k,\ell) \le (c/\ell)^{\ell/2} k^{k+3\ell/2-2}.
\end{equation}
Hence, if $k\le b n\qww$ for some fixed $b$, and $n$ is so large that
$p/(1-p)<2/n$,
\begin{equation}
  \label{pb}
\pkp(\ell) 
\le 
\frac{C(k,\ell)}{C(k,0)} \parfrac{p}{1-p}^\ell 
\le
\parfrac{4cb^3}{\ell}^{\ell/2}.
\end{equation}

Consider a sequence $k=k(n)$ such that $k n\qiww\to x$ for some $x>0$.
By \eqref{wright}, for every $\ell\ge0$, as \ntoo,
\begin{equation}
  \label{z5a}
\frac{C(k,\ell)}{C(k,0)} \parfrac{p}{1-p}^\ell 
=\wl k^{3\ell/2} n^{-\ell}\bigpar{1+o(1)}
\to\wl x^{3\ell/2}.
\end{equation}
Together with \eqref{pb}, this implies by dominated convergence
\begin{equation}
  \label{z5b}
\sum_{\ell=0}^\infty\frac{C(k,\ell)}{C(k,0)} \parfrac{p}{1-p}^\ell 
\to\sum_{\ell=0}^\infty\wl x^{3\ell/2}
=\Psi\bigpar{x^{3/2}}.
\end{equation}
Consequently, from \eqref{pkp}, \eqref{z5a},  \eqref{z5b} and
\eqref{pxl},
for \ntoo{} and every fixed $\ell$,
\begin{equation*}
  \pkp(\ell) \to   P_x(\ell).
\end{equation*}
Thus, the distribution $\pkp$ converges to $P_x$.

Let $\Xi'$ be the labelled point process constructed in the statement
of the theorem.
It follows from \refL{L:point2} and 
the Skorohod coupling theorem, see \eg{} \cite[Theorem 4.30]{Kallenberg},
that we may assume $\Xi_n$ and $\Xi$ to be coupled such that
$\xini\to\xi_i$ \as{} for every $i$.
By the description of $\Xix$ above and the convergence of $\pkp$ to
$P_x$ when $kn\qiww\to x$, it follows that we may couple also the
labels such that $\Xix_n\asto \Xi'$. 
Hence $\Xix_n\dto \Xi'$, and thus $\Xi'\eqd\Xix$.  
\end{proof}

\section{Intensity}\label{Sint}

Let, changing the notation slightly from \cite{LPW}, 
$X_n(k;\ell)$ denote the number of components with
$k$ vertices and complexity $\ell$ in the random graph \gnp, and let
$Y_n(k;\ell):=kX_n(k;\ell)$, the number of vertices in these components.
We further define,
for an interval $I$, 
$X_n(I;\ell):=\sum_{k\in I}X_n(k;\ell)$,
$X_n(k;\gexx\ell):=\sum_{j=\ell}^\infty X_n(k;j)$ 
and
$X_n(I;\gexx\ell):=\sum_{j=\ell}^\infty X_n(I;j)$,
and similary for $Y$. 
Thus, for example, $X_n(k)=\sum_j X_n(k;j)=X_n(k; \gex0)$.

Consider now a fixed $\ell\ge0$ and $k\le Cn^{2/3}$ for an arbitrary
constant $C$. 
Then, by well-known calculations, uniformly for all such $k$,
\begin{equation}
\label{sofie}
\begin{split}
\hskip1em&\hskip-1em
\E X_n(k,\ell)
=
\binom nk C(k,\ell) p^{k+\ell-1} (1-p)^{(n-k)k+\binom k2 -k-\ell+1}\\ 
&=\frac{n^k}{k!} \exp\Bigpar{-\frac{k^2}{2n}-\frac{k^3}{6n^2}+O\parfrac kn}
C(k,\ell)n^{1-k-\ell}(1+\gl n^{-1/3})^{k+\ell-1}
\\&\hskip22em \times
(1-p)^{nk-k^2/2-3k/2-\ell+1}\\
&=n^{1-\ell}\frac{C(k,\ell)}{k!} 
 \exp\Bigpar{-\frac{k^2}{2n}-\frac{k^3}{6n^2}
+(k+\ell-1)\gl n^{-1/3}-\frac12k\gl^2n^{-2/3}
\\&\hskip10em
-k-k\gl n^{-1/3}
+\frac{k^2}{2n}+\frac12k^2\gl n^{-4/3}
+O\parfrac kn +O(n^{-2/3})}\\
&=n^{1-\ell}\frac{C(k,\ell)}{k!} 
 \exp\Bigpar{-k-\frac{k^3}{6n^2}
-\frac12k\gl^2n^{-2/3}
+\frac12k^2\gl n^{-4/3}
+(\ell-1)\gl n^{-1/3}
\\&\hskip20em
+O\parfrac kn +O(n^{-2/3})}\\
&=n^{1-\ell}\frac{C(k,\ell)}{k!} 
e^{-k} 
 \exp\bigpar{-F(kn^{-2/3},\gl)}
\Bigpar{1+(\ell-1)\gl n^{-1/3}+O\parfrac kn +O(n^{-2/3})}
\end{split}
\end{equation}
where 
\begin{equation}\label{f}
F(x,\gl):=\frac16x^3-\frac12x^2\gl+\frac12x\gl^2=\frac{(x-\gl)^3+\gl^3}6.
\end{equation}
Note that
\begin{equation}
\label{magnus}
F(x,\gl)=\frac{x^3}{24}+x\frac{(x-2\gl)^2}{8}
\ge\frac{x^3}{24}
\ge0
\end{equation}
for all $x\ge0$ and $-\infty<\gl<\infty$.

In our first application of \eqref{sofie}, assume $0<a<b<\infty$ and
consider only $k\in[an\qww,bn\qww]$. For such $k$ and fixed $\ell$,
\eqref{sofie} gives, by \eqref{wright} and Stirling's formula, 
\begin{equation*}
  \begin{split}
  \E X_n(k,\ell)
&\sim
n^{1-\ell}\wl(2\pi)\qqi k^{3\ell/2-5/2}
e^{-F(kn^{-2/3},\gl)}
\\&
=
(2\pi)\qqi\wl\parfrac{k^{3/2}}{n}^{\ell-1}
e^{-F(kn^{-2/3},\gl)}k\qi,
  \end{split}
\end{equation*}
and summing over $k$ we obtain, as \ntoo,
\begin{equation}\label{erika}
  \E\bigpar{\xil_n[a,b]}=\E\sum_{k=an\qww}^{bn\qww}X_n(k,\ell)
\to
(2\pi)\qqi\wl
\int_a^b \bigpar{x^{3/2}}^{\ell-1}
e^{-F(x,\gl)}\,\frac{\dd x}x
.
\end{equation}

Since $\xil_n\dto\xil$ we have, by \refL{L:point3},
$\xil_n[a,b]\dto\xil[a,b]$ whenever $a$ and $b$ are continuity points
of $\xil$. In this case, by Fatou's lemma, $\E\xil[a,b]$ is at most
the right hand side of \eqref{erika}.

For any $a\in(0,\infty)$, $a\pm\eps$ are continuity points of $\xil$ 
for all but at most countably many $\eps\in(0,a)$, and for such
$\eps$ we thus obtain a formula for
$\E(\xil[a-\eps,a+\eps])$ and thus an upper bound of
$\E(\xil\set{a})$.
Letting $\eps\to0$ through such $\eps$, we see that 
$\E(\xil\set{a})=0$, so every point is a continuity point.
Consequently,
$\xil_n[a,b]\dto\xil[a,b]$ whenever $0<a<b\le\infty$.
Summing over all $\ell$, we see that every point is a continuity point
of $\Xi$ too, and thus
$\Xi_n[a,b]\dto\Xi[a,b]$ whenever $0<a<b\le\infty$.

To prove convergence of the expectations, we verify uniform
integrability by considering second moments.
(See also the more general \refL{Lmoments} below; we will give
a more elementary argument here, which in any case will be needed later.)

For simplicity, fix $\ell$ and write $\ek\=\E X_n(k,\ell)$.
Further, let $\ekj$ denote the expected number of ordered pairs of distinct
components of complexity $\ell$, of orders $k$ and $j$, respectively,
in \gnp.
Thus, if $k\neq j$ then $\ekj=\E\bigpar{X_n(k;\ell) X_n(j;\ell)}$,
while 
$$
\ekk=\E\bigpar{X_n(k;\ell)(X_n(k;\ell)-1)}=\E\bigpar{X_n(k;\ell)}^2-\ek.$$
Consequently, 
\begin{equation}\label{ex2}
\E\bigpar{ X_n\bigpar{\iab;\ell}}^2
=\E\biggpar{\sumkab X_n(k;\ell)}^2
=\sumkab\sumjab  \ekj+ \sumkab \ek.
\end{equation}
We have, cf.\ \eqref{sofie}, by simple calculations, assuming, say,
$k+j\le n/2$, 
\begin{equation}\label{ekj}
\begin{split}
\ekj&=\binom n{k+j}\binom{k+j}k C(k,\ell)C(j,\ell)
 p^{k+\ell-1+j+\ell-1}(1-p)^{n(k+j)-(k+j)^2/2-3(k+j)/2-2\ell+2}\\
&=\frac{\fall n{k+j}}{{\fall nk}{\fall nj}} (1-p)^{-kj} \ek\ej\\
&=\ek\ej\exp\Bigpar{\sum_{i=0}^{j-1}\ln\Bigpar{1-\frac k{n-i}}-kj\ln(1-p)}\\
&=\ek\ej\exp\Bigpar{\sum_{i=0}^{j-1}
-\Bigpar{\frac k{n}+\frac {ki}{n^2}+\frac {k^2}{2n^2}+
O\parfrac{ki^2+k^2i+k^3}{n^3}}
+kjp
 +O\parfrac{kj}{n^2}}\\
&=\ek\ej\exp\Bigpar{\gl{kj}n^{-4/3}-\frac{k^2j+kj^2}{2n^2}
 +O\Bigpar{\frac{kj}{n^2}+\frac{kj(k^2+j^2)}{n^3}}}.
\end{split}
\raisetag\baselineskip
\end{equation}

In particular, for $k,j\le C n\qww$, $\ekj=O(\ek\ej)$, and 
\eqref{ex2} implies, for fixed $\ell$, $a$ and $b$, with
$0<a<b<\infty$, 
recalling that $\E\xil_n[a,b]=O(1)$ by \eqref{erika},
\begin{equation*}
  \begin{split}
\E\bigpar{\xil_n[a,b]}^2
&
=
\E\bigpar{ X_n\bigpar{\iab;\ell}}^2
\\&
=
O\Bigpar{\bigpar{\E X_n\bigpar{\iab;\ell}}^2
+ \E X_n\bigpar{\iab;\ell}}  
\\&
=
O\Bigpar{\bigpar{\E\xil_n[a,b]}^2
+ \E \xil_n[a,b]}  
=O(1).	
  \end{split}
\end{equation*}
Thus, the random variables $\xil_n[a,b]$ are uniformly integrable, and
$\xil_n[a,b]\dto\xil[a,b]$ implies 
$\E\xil_n[a,b]\to\E\xil[a,b]$, 
see \eg{} \cite[Theorems 5.4.2 and 5.5.9]{Gut}.
Consequently, 
$\E\xil[a,b]$ equals the \rhs{} of
\eqref{erika}. This leads to the following result.
Recall that $\Psi$ denotes the moment generating function 
\eqref{mgf}
of the Brownian excursion area.

\begin{theorem}\label{Tint}
The point process
$\xil$ has intensity 
$\gL_\ell\=(2\pi)\qqi\wl x^{3\ell/2-5/2} e^{-F(x,\gl)}$
on $(0,\infty)$.
Their sum $\Xi$ has the intensity, for $0<x<\infty$,
\begin{equation}\label{tint}
\gL(x)=
\gL\gll(x)\=
\sum_{\ell=0}^\infty\gL_\ell(x)
=
(2\pi)\qqi x^{-5/2} \Psi(x^{3/2}) e^{-F(x,\gl)}.
\end{equation}
\end{theorem}

\begin{proof}
We have shown that
$\E\xil[a,b]=\int_a^b \gL_\ell(x)\dd x$ when $0<a<b<\infty$, which
by definition shows that
$\gL_\ell(x)$ is the intensity of $\xil$. The second part follows by
summing over $\ell$.
\end{proof}

\begin{corollary}
  \label{C1}
\begin{equation}\label{c1}
\E\ze
=\int_\eps^\infty
(2\pi)\qqi x^{-3/2} \Psi(x^{3/2}) e^{-F(x,\gl)}\dd x.
\end{equation}
\end{corollary}

\begin{proof}
$  \E\ze=\int_\eps^\infty x \gL(x)\dd x$.
\end{proof}

We already know that the expectation in \eqref{c1} is finite; that the integral
converges follows also by
\eqref{jesper} and \eqref{magnus}, which imply that
$\gL(x)$ decreases exponentially as \xtoo{}.
Note further that the intensity 
$\gL\gll(x)\sim (2\pi)\qqi x^{-5/2}$ as $x\to0$, for every $\gl$.

\begin{remark}\label{Rx}
The intensity $\gL_\ell$ has a finite integral over $(0,\infty)$
precisely when the exponent $3\ell/2-5/2> -1$. 
Thus, for $\ell=0$ and $\ell=1$, $\xil$ has an infinite expected 
number of points; indeed, it is easily shown from \eqref{vx0} and
\eqref{vx1} that $\xil$ \as{} has an infinite number of points.
On the other hand,
for any $\ell\geq 2$, $\xil$ has a finite number of points. Further, 
$\sum_{\ell\geq 2}\xil$, the point process for the complex components, 
has a finite number of points.  One may view this in an evolutionary way.  
Roughly speaking, when $\lambda$ is large negative complex components 
have not yet formed.  When $\lambda$ is large positive a ``dominant
component" will have formed which is complex.  But there will not
usually be other complex components as components get
``sucked into" the dominant component before becoming complex.
In \cite{giant} it is shown that with probability converging to
$5\pi/18\approx 0.87$
there is never more than one complex component  in the entire
evolution of the random graph.
\end{remark}

\begin{remark}
Considering the difference $\E\ze\gll-\E\ze^{(0)}$ and   letting
$\eps\to0$, we find from \refT{T1} and \refC{C1} the identity
\begin{equation}
  \label{sjw}
\intoo
(2\pi)\qqi x^{-3/2} \Psi(x^{3/2})\Bigpar{e^{-F(x,\gl)}-e^{-F(x,0)}}\dd x
=\intoo x\bigpar{\gL\gll(x)-\gL\gxx0(x)}\dd x
=\gl.
\end{equation}
Differentiating with respect to $\gl$ we further find
$-\intoo x\frac{\partial F}{\partial\gl}(x,\gl)\gL\gll(x)\dd x=1$
or
\begin{equation*}
  \intoo x^2(x-2\gl)\gL\gll(x)\dd x=2,
\end{equation*}
and thus $\E\sum_i\xi_i^3=2+2\gl\E\sum_i\xi_i^2$.
\end{remark}

\begin{remark}
We similarly find expressions for the expectations of sums of all
points in $\Xix$ with a given label. For example, for label 1,
corresponding to unicyclic components, we obtain 
for the expectation of the total weight of $\xilx1$
\begin{equation}
\label{manne}  
\intoo x\gL_1(x)\dd x=\tfrac14\intoo e^{-F(x,\gl)}\dd x.
\end{equation}
Thus, \cf{} \refR{Rx},
the $\xilx1$ process has an infinite number of points with finite sum. 
See also \cite[Lemma 2.2]{LPW}, which implies both \eqref{manne} and
\begin{equation*}
  \intoo x\bigpar{\gL\gll(x)-\gL\gll_0(x)}\dd x
=(2\pi)\qqi\intoo x^{-3/2}\bigpar{1-e^{-F(x,\gl)}}\dd x+\gl,
\end{equation*}
which indeed also easily follows from \eqref{t1a} and \eqref{c1}.
The expectation of the sum of all points with label at least 2
(corresponding to the total size of the complex components in \gnp)
is  $\intoo x\bigpar{\gL\gll(x)-\gL\gll_0(x)-\gL\gll_1(x)}\dd x$; 
an evaluation in terms of hypergeometric functions is given in
\cite[(15.12)]{giant}. 
\end{remark}

\section{An estimate for \gnp} \label{Smoments} 

We prove in this section an estimate for the components of \gnp{} that
we will need. This estimate is known, at least in principle, but we do
not know any precise reference.

\begin{lemma}
  \label{Lmoments}
For any fixed $\eps>0$ and $q\ge0$,
\begin{align*}
\E\bigpar{\bigpar{X_n[\eps n\qww,\infty)}^q}
&=O\bigpar{1},	
&
\E\bigpar{\bigpar{Y_n[\eps n\qww,\infty)}^q}
&=O\bigpar{n^{2q/3}}.	 
\end{align*}
In other words, $\wne$, $\zne$  and $\znex$ have moments that are
bounded, uniformly in $n$. 
\end{lemma}

\begin{proof}
It suffices to prove the result for $Y_n$, since
$Y_n[\eps n\qww,\infty) \ge \eps n\qww X_n[\eps n\qww,\infty)$.

Let us begin with the complex components; in this case we do not
need a lower bound on the size of the components. (This is not
surprising, since typically there are no small complex components.)
Let $\nc(G)$ denote the number of vertices in complex components of a
graph $G$.
Thus $\nc\bigpar{\gnp}=Y_n\bigpar{[1,\infty),\gex2}$.

The \emph{excess} of a graph, as defined in \cite[\S13]{giant}, equals
the complexity minus the number of complex components. (Thus, a
component of complexity $\ell$ contributes $\max(\ell-1,0)$.)
We first claim that there exists $\eta>0$ such that if $q\ge1$ and
$\mu\in\bbR$ are fixed, and we momentarily consider the random graph
\gnm{} with a fixed number of edges $m=\floor{\frac n2(1+\mu n\qiw)}$,
then
\begin{equation}
  \label{r1}
\E\Bigpar{\nc\bigpar{\gnm}^q
 \, \ett\bigl[\operatorname{excess}\bigpar{\gnm}=r\bigr]}
=O\bigpar{n^{2q/3}(r+1)^qe^{-\eta r}},
\end{equation}
uniformly in $n$ and $r\ge0$.

Indeed, the case $q=0$ of \eqref{r1} is a special case of 
\cite[Lemma 5]{giant} (with $d=0$), and the general case follows by a simple
modification of the (not so simple) proof, as is remarked for the
case $q=1$ on \cite[pages 299--300]{giant}.

We may thus sum \eqref{r1} over $r\ge0$, and find, still for fixed $q$
and $\mu$,
\begin{equation}
  \label{r2}
\E\bigpar{\nc\bigpar{\gnm}^q}
=O\bigpar{n^{2q/3}}.
\end{equation}
Given $\gl$, choose $\mu>\gl$ and observe that by a standard Chernoff estimate,
see for example \cite[Theorem 2.1]{JLR}, the probability that \gnp{}
has more than $m$ edges is $O(e^{-\gd n\qw})$ for some $\gd>0$.
Since $\nc$ is monotone if we add edges, any coupling of $\gnp$ and
$\gnm$ thus gives, for fixed $q\ge0$,
\begin{equation}
  \label{r3}
  \begin{split}
\E\Bigpar{Y_n\bigpar{[1,\infty),\gex2}^q}
&
=
\E\bigpar{\nc\bigpar{\gnp}^q}
\le
\E\bigpar{\nc\bigpar{\gnm}^q}+O\bigpar{n^q e^{-\gd n\qw}}
\\&
=O\bigpar{n^{2q/3}}.	
  \end{split}
\end{equation}

For components of complexity 0 or 1, \ie{} trees and unicyclic
components, it is possible to argue as for the second moment in
\refS{Sint}, but we will instead use a trick together with the result
just proved.

Let $\pkk{\ell}$ be the probability that $j$ given disjoint subsets of
the vertex set of \gnp, with sizes $k_1,\dots,k_j$ respectively,
all are the vertex sets of components with complexities $\ell$.
Thus, with $k=k_1+\dots+k_j$,
\begin{equation}\label{pkk}
  \pkk{\ell}
=
(1-p)^{nk-k^2/2-3k/2-j\ell+j}
\prod_{i=1}^j {C(k_i,\ell) p^{k_i+\ell-1}}.
\end{equation}
It is easily seen that for any integer $q\ge1$,
\begin{equation}\label{pku}
  \E\bigpar{Y_n\bigpar{[A,\infty);\ell}^q}
=
\sum_{j=1}^q\sum_{k_1,\dots,k_j\ge A} c(n;q;j;k_1,\dots,k_j)\pkk{\ell},
\end{equation}
for some combinatorial coefficients $ c(n;q;j;k_1,\dots,k_j)$ not
depending on $\ell$.
For fixed $\ell$ and $\ell'$, we have by \eqref{pkk} and \eqref{wright},
\begin{equation*}
  \frac{\pkk{\ell'}}{\pkk{\ell}}
=
\prod_{i=1}^j \frac{C(k_i,\ell')}{C(k_i,\ell)}
\parfrac{p}{1-p}^{\ell'-\ell}
=\Theta{\biggpar{\prod_{i=1}^j \frac{k_i^{3(\ell'-\ell)/2}}{n^{\ell'-\ell}}}}.
\end{equation*}
Hence, if $\ell'\le\ell$ and $k_i\ge \eps n\qww$, we have
$\pkk{\ell'}=O\bigpar{\pkk{\ell}}$
(recall that $\eps$ is fixed), and \eqref{pku} yields
\begin{equation*}
    \E\bigpar{Y_n\bigpar{[\eps n\qww,\infty);\ell'}^q}
=
O\Bigpar{
  \E\bigpar{Y_n\bigpar{[\eps n\qww,\infty);\ell}^q}
}
\end{equation*}
We apply this with $\ell'=0$ and $1$ and $\ell=2$, and obtain
from \eqref{r3} the required estimates for 
$\E\bigpar{Y_n\bigpar{[\eps n\qww,\infty);0}^q}$ and
$\E\bigpar{Y_n\bigpar{[\eps n\qww,\infty);1}^q}$, which together with
	\eqref{r3} complete the proof.
\end{proof}

Before we proceed, we point out a simple consequence.
Peres \cite{Peres} recently gave a simple proof (with an explicit
bound) of the case $q=2$; this case is equivalent to 
$\E |\cC(v)|=O\bigpar{n^{1/3}}$ where $\cC(v)$ is the component containing
a given (or random) vertex $v$.

\begin{corollary}
  \label{Cp}
Let $q>3/2$. Then $\E\sum_i\xini^q=O(1)$;
equivalently, for \gnp,
$\E\sum_i|\cC_i|^q=O\bigpar{n^{2q/3}}$ for any $q>3/2$.
\end{corollary}

\begin{proof}
First,
$\sum_{i:|\cC_i|>n\qww}|\cC_i|^q
\le Y_n[n\qww,\infty)^q,$
whose mean is $O\bigpar{n^{2q/3}}$ by \refL{Lmoments}.

Similarly, the sum over the complex components has expectation
$O\bigpar{n^{2q/3}}$ by \eqref{r3}.
It thus remains only to consider components of size at most $n\qww$
with complexity 0 or 1.
The corresponding sum has expectation
$
  \sum_{k=1}^{n\qww} k^q(t_k+u_k)
$, 
which is $O\bigpar{n^{2q/3}}$ by \eqref{tk1} and \eqref{ut}.
\end{proof}

\begin{proof}[Proof of \refP{P1}]  
It is an easy consequence of
\refL{L:point2} 
that the mappings $\set{\xi_i}_i\mapsto \sum_{\xi_i\ge\eps}\xi_i$,
$\set{\xi_i}_i\mapsto \sum_{\xi_i>\eps}\xi_i$ and
$\set{\xi_i}_i\mapsto \#\set{i:\xi_i\ge\eps}$
are measurable on $\fN(0,\infty]$ and continuous at every 
$\set{\xi_i}$ such that $\eps\notin\set{\xi_i}$.

Since $\Xi_n\dto\Xi$ and $\P(\eps\in\Xi)=0$ by \refT{Tint},
the results
$\zne\dto\ze$, $\znex\dto\ze$ and $\wne\dto\we$ follow
by the continuous mapping theorem, see \eg{} \cite[Theorem 5.1]{Billingsley}.

Since the moments of $\zne$, $\znex$ and $\wne$ are bounded uniformly
in $n$ by \refL{Lmoments}, this further implies convergence of all
moments,
see \eg{} \cite[Theorems 5.4.2 and 5.5.9]{Gut}.
\end{proof}

\section{Second partial proof of \refT{T1}} \label{Spf2}
In this proof we do the calculations with the small components, and 
consider complexities 0 and 1 separately.
Let throughout $0<\eps<1$.

Consider first the tree components.
Let $t_k=\E X_n(k;0)$ be the expected number of tree components of
order $k$.
By \eqref{sofie}, for
$k\le n^{2/3}$,
\begin{equation}
\label{tk}
\begin{split}
t_k
&=n\frac{k^{k-2}}{k!}e^{-k} 
 \exp\bigpar{-F(kn^{-2/3},\gl)}
\Bigpar{1-\gl n^{-1/3}+O\parfrac kn +O(n^{-2/3})}.
\end{split}
\end{equation}

In particular, with $\tkx\=n{k^{k-2}}e^{-k}/{k!}$,
\begin{equation}\label{tk1}
t_k
=
O(\tkx)
=
O\Bigpar{n\frac{k^{k-2}}{k!}e^{-k}}
=O\parfrac{n}{k^{5/2}}. 
\end{equation}
Note further that, 
for any fixed real $\ga$ and all $\eps>0$,
\begin{equation}
\label{tk2}
\sumke k^\ga \tkx =
\sumke O\bigpar{n k^{\ga-5/2}} =
\begin{cases}
O\bigpar{n(\eps n^{2/3})^{\ga-3/2}} 
= O\bigpar{n^{2\ga/3} \eps^{\ga-3/2}},
&\ga>3/2,\\
O(n),&\ga<3/2.
\end{cases}
\end{equation}

By \eqref{tk} and \eqref{tk2} we obtain, 
\begin{equation}
\label{tk3}
\begin{split}
\E Y([1,\eps n\qww],0)
&=
\sumke kt_k
=n\sumke \frac{k^{k-1}}{k!}e^{-k} 
\bigpar{1-\gl n^{-1/3}}
\\&\hskip8em
+O\Bigpar{\sumke k\tkx\bigpar{F(kn\qiww,\gl)+kn^{-1} +n^{-2/3}}}\\
&=\bigpar{n-\gl n^{2/3}}\sumke \frac{k^{k-1}}{k!}e^{-k} 
+O\Bigpar{\sumke \tkx(k^4n^{-2} 
 +k^2n^{-2/3})}\\
&=\bigpar{n-\gl n^{2/3}}\sumke \frac{k^{k-1}}{k!}e^{-k} 
+O\Bigpar{\eps^{1/2}n^{2/3}}.
\raisetag{24pt} 
\end{split}
\end{equation}
Moreover, using the fact that 
$\sum_1^\infty \frac{k^{k-1}}{k!}e^{-k} =1$,
and Stirling's formula,
\begin{equation}
\label{tk4}
\begin{split}
\sumke \frac{k^{k-1}}{k!}e^{-k} 
&=1-\sum_{k>\eps n^{2/3}} \frac{k^{k-1}}{k!}e^{-k}\\
&=1-\sum_{k>\eps n^{2/3}} (2\pi k^3)^{-1/2}
+O\Bigpar{\sum_{k>\eps n^{2/3}} k^{-5/2}}\\
&=1-\int_{\eps n^{2/3}}^\infty (2\pi x^3)^{-1/2}\,dx
+O(\eps^{-3/2} n^{-1})\\
&=1-\sqrt{\tfrac2{\pi}}\eps^{-1/2}n^{-1/3}
+O(\eps^{-3/2} n^{-1}).
\end{split}
\end{equation}
Consequently, combining \eqref{tk3} and \eqref{tk4},
\begin{multline}\label{e0}
\E Y_n\bigpar{[1,\eps n^{2/3}];0}
\\
=n-\sqrt{\tfrac2{\pi}}\eps^{-1/2}n^{2/3}-\gl n^{2/3}
+O\Bigpar{\eps^{1/2}n^{2/3}+\eps^{-3/2}+\eps^{-1/2}n^{1/3}}.
\end{multline}

Next,
let $u_k=\E X_n(k;1)$ be the expected number of unicyclic components of
order $k$.
We have, \cf{} \eqref{sofie} and \eqref{wright},
\begin{equation}\label{ut}
\begin{split}
u_k&=\binom nk C(k,1) p^k(1-p)^{(n-k)k+\binom k2-k}\\
&=\frac{C(k,1)}{k^{k-2}}p(1-p)^{-1} t_k\\
&=O\Bigpar{n^{-1}k^{3/2}t_k},
\end{split}
\end{equation}
and thus, by \eqref{tk2},
\begin{equation}\label{e1}
\E Y_n\bigpar{[1,\eps n^{2/3}];1}=
\sumke ku_k
= O\Bigpar{n^{-1}\sumke k^{5/2}t_k} = O(n^{2/3}\eps).
\end{equation}

For complex components we use the well-known fact that
$\E X_n\bigpar{[1,\infty);\ge\nobreak2}$ is bounded; 
see the stronger result in  
\cite{SJ:multi}, \cite[Theorem 5.8(i)]{JLR}. 
(As a bound we can take 1.2, say, at least for large
$n$, and possibly 1, as conjectured in \cite{LPW}.)
Hence,
\begin{equation}\label{e2}
\E Y_n\bigpar{[1,\eps n^{2/3}];\gex2}
\le
\eps n^{2/3} \E X_n\bigpar{[1,\eps n^{2/3}];\gex2}
= O(n^{2/3}\eps).
\end{equation}

Adding  \eqref{e0}, \eqref{e1} and \eqref{e2}, we find, since the sum
of all component sizes $Y_n([1,\infty))=n$,
\begin{multline*}
\E\znex= 
n\qiww\E Y_n\bigpar{\eps n^{2/3},\infty}
=
n\qiww\bigpar{n-\E Y_n\bigpar{[1,\eps n^{2/3}]}}
\\
=\sqrt{\tfrac2{\pi}}\eps^{-1/2}+\gl 
+O\Bigpar{\eps^{1/2}+\eps^{-3/2}n\qiww+\eps^{-1/2}n\qiw}.
\end{multline*}
Thus, letting \ntoo, by \refP{P1},
\begin{equation*}
  \E\ze
=\sqrt{\tfrac2{\pi}}\eps^{-1/2}+\gl 
+O\Bigpar{\eps^{1/2}},
\end{equation*}
which is \eqref{t1a} with the weaker error term $O(\eps\qq)$.

Next, consider the variance of $Y_n\bigpar{[1,\eps n^{2/3}];0}$.
Similarly to \eqref{ex2} we have, with $\ell=0$, and $\ek=t_k$,
\begin{equation}\label{ey2}
\E\bigpar{ Y_n\bigpar{\iie;\ell}}^2
=\E\biggpar{\sumke kX_n(k;\ell)}^2
=\sumke\sumje kj \ekj+ \sumke k^2\ek.
\end{equation}
Hence,
using \eqref{ekj} 
and letting $A_2\=\sumke k^2 t_k=O\bigpar{ n^{4/3}\eps\qq}$, 
by \eqref{tk2}, 
\begin{equation}\label{yvar0}
\begin{split}
\Var\bigpar{ Y_n\bigpar{[1,\eps n^{2/3}];0}}
&
= 
\sumke k^2t_k+\sumke\sumje kj (\ekj-t_kt_j)
\\&
=
\sumke k^2t_k+\sumke\sumje kj t_kt_j 
O\Bigpar{{kj}n^{-4/3}}\\
&= \sumke k^2t_k+
O\Bigpar{\sumke\sumje k^2j^2 t_kt_j n^{-4/3}}\\
&=A_2+O\bigpar{A_2^2n^{-4/3}}
=A_2+O\bigpar{ n^{4/3}\eps}.
\end{split}
\end{equation}
In particular, by \eqref{tk2}, this variance is $O\bigpar{ n^{4/3}\eps\qq}$.

The variance of $Y_n\bigpar{[1,\eps n^{2/3}];1}$ can be computed in the
same way, with $t_k$ replaced by $u_k$. Since $u_k=O(\eps^{3/2}t_k)$ for
$k\le\eps n\qww$ by \eqref{ut}, we obtain the estimate
\begin{equation}\label{yvar1}
\Var\bigpar{ Y_n\bigpar{[1,\eps n^{2/3}];1}}
= 
O\bigpar{\eps^{3/2}A_2+\eps^3A_2^2n^{-4/3}}
=O\bigpar{ n^{4/3}\eps^2}.
\end{equation}

For the complex components we now use the fact that also
$\E\bigpar{ X_n\bigpar{[1,\infty);\gex2}^2}$ is bounded \cite{SJ:multi}, 
and thus
\begin{equation}\label{yvar2}
  \begin{split}
\Var\bigpar{ Y_n\bigpar{[1,\eps n^{2/3}];\gex2}}
&\le
\E\bigpar{ Y_n\bigpar{[1,\eps n^{2/3}];\gex2}^2}
\le
\eps^2 n^{4/3} \E \bigpar{ X_n\bigpar{[1,\eps n^{2/3}];\gex2}^2}
\\
&= O(n^{4/3}\eps^2).
  \end{split}
\raisetag\baselineskip
\end{equation}

By the \CS{} inequality, the three covariances between the three
variables in \eqref{yvar0}, \eqref{yvar1} and \eqref{yvar2} are all
$O(n^{4/3}\eps^{5/4})$, so summing the variables we find from these
formulas that 
\begin{equation}\label{yvar}  
\Var\bigpar{ Y_n(\eps n^{2/3},\infty)}
=
\Var\bigpar{ Y_n[1,\eps n^{2/3}]}
=
A_2+O\bigpar{ n^{4/3}\eps}.
\end{equation}
Moreover, by \eqref{tk}, \eqref{tk2} and Stirling's formula,
\begin{equation*}
  \begin{split}
A_2
&=
\sumke k^2t_k
=
\sumke k^2\tkx\bigpar{1+O(kn\qiww+n\qiw)}
\\&
=
\sumke n\bigpar{(2\pi k)\qqi+O(k^{-3/2})}
 +O\bigpar{n^{4/3}\eps^{3/2}+n\eps\qq}
\\&
=
 (2/\pi)\qq\eps\qq n^{1+1/3}
 +O\bigpar{n+n^{4/3}\eps^{3/2}}.
  \end{split}
\end{equation*}
Thus, \eqref{yvar} yields
\begin{equation*}  
\Var\znex= 
n^{-4/3}\Var\bigpar{ Y_n(\eps n^{2/3},\infty)}
=
(2/\pi)\qq\eps\qq +O\bigpar{\eps+n\qiw}.
\end{equation*}
and \eqref{t1b}
follows by \refP{P1}.


\section{Proof of \refT{T2}} 
By \refT{Tint} and \eqref{f}, for $0<\eps\le1$ and fixed $\gl$,
\begin{equation*}
  \begin{split}
\E\we
&=
\int_\eps^\infty \gL(x)\dd x
=\E W_1 +
 \int_\eps^1 (2\pi)\qqi x^{-5/2} \Psi(x^{3/2}) e^{-F(x,\gl)}\dd x
\\&
=
O(1)+
 \int_\eps^1 (2\pi)\qqi x^{-5/2}\bigpar{1+w_1x^{3/2}+O(x^3)} 
 \bigpar{1-F(x,\gl)+O(x^2)}\dd x
\\&
= \int_\eps^1 (2\pi)\qqi x^{-5/2}\bigpar{1+w_1x^{3/2}-x\gl^2/2+O(x^2)}\dd x
+O(1),
  \end{split}
\end{equation*}
and \eqref{t2a} follows. (Recall that $w_1=(\pi/8)\qq$.)

For the variance, we use \eqref{ex2} with $a=\eps$ and $b=1$ and argue
as in \refS{Spf2}.
Calculations similar to \eqref{yvar0} yield
\begin{align}\label{vx0}
\Var\bigpar{ X_n\bigpar{[\eps n^{2/3},n\qww];0}}
&= 
\E\bigpar{ X_n\bigpar{[\eps n^{2/3},n\qww];0}}+O(\eps\qi)
=O(\eps^{-3/2})
\intertext{and}
\label{vx1}
\Var\bigpar{ X_n\bigpar{[\eps n^{2/3},n\qww];1}}
&= 
\E\bigpar{ X_n\bigpar{[\eps n^{2/3},n\qww];1}}+O(1)
=O\bigpar{\ln(1/\eps)};
\end{align}
we omit the details. Using again the fact that
$\E\bigpar{ X_n\bigpar{[1,\infty];\gex2}^2}=O(1)$
together with
$\E\bigpar{\bigpar{X_n(n\qww,\infty]}^2}=O(1)$ 
(\refL{Lmoments})
and the \CS{} inequality,
we obtain
\begin{equation*}
\Var\bigpar{ X_n[\eps n^{2/3},\infty]}
= 
\E\bigpar{ X_n[\eps n^{2/3},\infty]}
+O(\eps\qi).
\end{equation*}
Letting \ntoo, we find using \refP{P1},
$\Var\we=\E\we+O(\eps\qi)$, and \eqref{t2b} follows.

\section{The Palm distribution} \label{Spalm}

The Palm distributions of a point process $\Xi$ in a suitable space
$\fS$ are the conditional distributions $\cL(\Xi\mid s\in\Xi)$
given the presence of a given point $s$, $s\in\fS$. (Usually,
$s\in\fS$ is an event of probability 0, so this must be interpreted
with some care, see \cite[Chapter 10]{Kallenberg0}.
In particular, note that the Palm distribution is uniquely determined
only for \aex{} $s$.)

In our case, the Palm distribution is obtained by a simple shift of
the parameter $\gl$; we thus write $\Xila=\set{\xi_i\gll}_i$ in this section.
Recall that we regard $\Xila$ as a random measure on $(0,\infty)$ that
is the sum of the pointmasses $\gd_{\xi_i\gll}$, see \refA{S:point}.

\begin{theorem}\label{Tpalm}
  The Palm distribution $\cL(\Xila\mid s\in\Xila)$ equals for every
  $s>0$ the distribution of $\Xilas+\gd_s$.
\end{theorem}

\begin{proof}
  Given that \gnp{} has a component of size $m$ on a certain set of
  vertices, the remainder of the graph is distributed as $G(n-m,p)$.
Hence, if $\fN=\fN(0,\infty]$ is 
the space of locally fimite integer-valued measures 
on $\fS=(0,\infty]$
defined in \refA{S:point},
and
$f:(0,\infty]\to\bbR$ and $g:\fN\to\bbR$ are bounded continuous functions 
and $f$ has compact support, then
\begin{equation}
  \label{palm1}
\E\Bigpar{g(\Xi_n)\intoo f\dd\Xi_n}
=
\E\Bigpar{g(\Xi_n)\sum_i f(\xini)}
=
\E\Bigpar{\sum_i h_n(\xini)},
\end{equation}
where 
$h_n(s)=f(s)\E g(\Xi_{n-sn\qww,p}+\gd_s)$.
If \ntoo{} and $s_n\to s$, then
$(n-s_nn\qww)p=1+(\gl-s+o(1))n\qiw$ and 
$\Xi_{n-s_n n\qww,p}\dto \Xilas$, and thus 
$h_n(s_n)\to h(s)\=f(s)\E g(\Xilas+\gd_s)$.
It follows, using \refL{L:point2} and \cite[Theorem 5.5]{Billingsley}, that 
$\sum_i h_n(\xini)\dto \sum_i h(\xi_i)$, and thus by 
\eqref{palm1} and dominated convergence, 
\begin{multline}\label{dadel}
\E\Bigpar{g(\Xila)\intoo f\dd\Xila}
=\limn
\E\Bigpar{g(\Xi_n)\intoo f\dd\Xi_n}
=
\limn
\E\Bigpar{\sum_i h_n(\xini)}
\\
=
\E\Bigpar{\sum_i h(\xi_i)}
=
\intoo h(s)\dd\E\Xila(s)
=
\intoo f(s)\E g(\Xilas+\gd_s)\dd\E\Xila(s),
\end{multline}
where $\dd\E\Xila(s)=\gL\gll(s)\dd s$ by \refT{Tint}.
It follows by a monotone class argument 
(\eg{} \cite[Theorem A.1]{SJIII})
that the first and last terms are equal for any bounded measurable
$g$, and the result follows, see
\cite[(10.2)]{Kallenberg0}.
\end{proof}

Note that Theorems \refand{Tpalm}{T1} imply that for any fixed $\gl$ and $s>0$,
for small $\eps$ (so that $\eps<s$),
$\E(\ze\gll\mid s\in \Xila)=\E\ze\glls+s=\E\ze\gll+O(\eps\qq)$.
Hence the existence of a certain point in $\Xi$ asymptotically 
does not influence
$\E\ze$ for small $\eps$, showing the rigidity of $\Xi$.

\refT{Tpalm} can be put in a computational form as follows. Let, as
above, $\fN=\fN(0,\infty]$ be the space of integer-valued measures 
defined in \refA{S:point}.

\begin{theorem}
  \label{Tpalm2}
For any bounded or non-negative
measurable function $F:(0,\infty)\times\fN\to[0,\infty]$,
\begin{equation}
  \label{palm}
\E\sum_i F\bigpar{\xi_i\gll,\Xila}
=\intoo\E F\bigpar{x,\Xi\gllx+\gd_x} \gL\gll(x)\dd x,
\end{equation}
where $\gL\gll(x)$ is given by \eqref{tint}.
\end{theorem}
\begin{proof}
First consider $F$ of the special form $F(x,\Xi)=f(x)g(\Xi)$, where, as
in the proof of \refT{Tpalm},
$f:(0,\infty]\to\bbR$ and $g:\fN\to\bbR$ are bounded continuous functions 
and $f$ has compact support. Then \eqref{dadel} holds, which can be
written
\begin{equation}\label{kokos}
  \E\intoo F(x,\Xila)\dd\Xila(x)
=
\intoo\E F(x,\Xi\gllx+\gd_x)\gL\gll(x)\dd x.
\end{equation}
By another monotone class argument 
(\eg{} \cite[Theorem A.1]{SJIII}), \eqref{kokos} holds for every bounded
measurable $F$, and thus by monotone convergence for every non-negative
measurable $F$ too.

The integral on the \lhs{} of \eqref{kokos} equals $\sum_i
F(\xi_i\gll,\Xila)$,
which yields \eqref{palm}.
\end{proof}

We give some applications.

\begin{corollary}
  \label{Cvar}
Let $\gL\gll(x)$ be given by \eqref{tint}. Then,
for every $\eps>0$,
\begin{align*}
  \E\ze^2 &
=\inteo x^2\gL\gll(x)\dd x
+\inteo \inteo xy\gL\gll(x)\gL\gllx(y)\dd y\dd x
\\
\intertext{and thus}
  \Var\ze &
=\inteo x^2\gL\gll(x)\dd x
-\inteo \inteo xy\gL\gll(x)\bigpar{\gL\gll(y)-\gL\gllx(y)}\dd y\dd x
\\&
=\inteo \intoe xy\gL\gll(x)\bigpar{\gL\gll(y)-\gL\gllx(y)}\dd y\dd x.
\end{align*}
\end{corollary}
\begin{proof}
Take $F(x,\Xi)=x\inteo y\dd \Xi(y)\ett[x\ge\eps]$ in \eqref{palm},
or $f(x)=x\ett[x\ge\eps]$ and 
$g(\Xi)=\inteo y\dd \Xi(y)$ in \eqref{dadel},
to find
\begin{equation*}
  \E\ze^2
=
\inteo x\E\Bigpar{\inteo y\dd\Xi\gllx(y)+x}\gL\gll(x)\dd x,
\end{equation*}
which yields the formula for $\E\ze^2$ by \refT{Tint} (or \refC{C1}) applied
with $\gl-x$.

The first formula for $\Var\ze$ follows immediately, and the second
  follows because \eqref{sjw} implies
  \begin{equation*}
x=\intoo y\bigpar{\gL\gll(y)-\gL\gllx(y)}\dd y	
  \end{equation*}
and thus
  \begin{equation*}
\inteo x^2\gL\gll(x)\dd x
=\inteo\intoo xy\bigpar{\gL\gll(y)-\gL\gllx(y)}\gL\gll(x)\dd y \dd x.
  \end{equation*}
\vskip-\baselineskip
\end{proof}

\begin{corollary}
  \label{Csimple}
$\Xila$ is \as{} simple, \ie{} lacks multiple points.
\end{corollary}

\begin{proof}
  Take $F(x,\Xi)\=\ett[\Xi\set{x}\ge2]$ in \eqref{palm}.
The \lhs{} becomes the expected number of multiple points (with
multiplicities), while the \rhs{} is 0 because, for each $x$, 
$\E F(x,\Xi\gllx+\gd_x)=\P\bigpar{\Xi\gllx\set{x}\ge1}=0$,
using \refT{Tint} which shows that the intensity of $\Xi\gllx$ is
absolutely continuous.
\end{proof}

\begin{corollary}
  \label{Cxi1}
The largest point $\xi_1\gll$ in $\Xila$ has a distribution with 
the density function $\hla1(x)\=\P\bigpar{\Xi\gllx(x,\infty)=0}\gL\gll(x)$.
\end{corollary}

\begin{proof}
Let $f:(0,\infty)\to[0,\infty]$ be a measurable function and 
  take $F(x,\Xi)\=f(x)\ett[\Xi(x,\infty)=0]$ in \refT{Tpalm2}.
Since $\Xila$ is simple by \refC{Csimple}, the \lhs{} of \eqref{palm}
becomes $\E f(\xi_1\gll)$, and the \rhs{} is
$\intoo f(x)\hla1(x)\dd x$.
Since $f$ is arbitrary, the result follows.
\end{proof}

The proof immediately extends to the following, more general, result.
\begin{corollary}
  \label{Cxik}
For any $k\ge1$, 
the $k$:th largest point $\xi_k\gll$ in $\Xila$ has a distribution with 
the density function
$\hla{k}(x)\=\P\bigpar{\Xi\gllx(x,\infty)=k-1}\gL\gll(x)$.
\nopf
\end{corollary}

\begin{corollary}
  \label{Cmom}
For any Borel set $\cA\subseteq(0,\infty)$ and $k\ge1$,
\begin{equation*}
  \E\bigpar{\fall{\Xila(\cA)} k} 
=
\int_\cA\dotsi\int_\cA \gL\gll(x_1)\gL\gxx{\gl-x_1}(x_2)
\dotsm
\gL\gxx{\gl-x_1-\dots-x_{k-1}}(x_k)\dd x_k\dotsm \dd x_1.
\end{equation*}
\end{corollary}

\begin{proof}
For $k=1$, this is just the definition of intensity, see
\refT{Tint}.
For $k\ge2$, we use \refT{Tpalm2} with 
$F(x,\Xi)\=\ett[x\in \cA]\fall{\bigpar{\Xi(\cA)-1}}{k-1}$,
which yields
\begin{equation*}
  \E\bigpar{\fall{\Xila(\cA)} k}
=\intoo\ett[x\in \cA]\E\fall{\bigpar{\Xi\gllx(\cA)}}{k-1}\gL\gll(x)\dd x,
\end{equation*}
and the result follows by induction.
\end{proof}

\begin{remark}
  \label{Rmom}
It follows immediately that if $\cA\subseteq[a,b]$ with $0<a<b<\infty$,
then $\E\bigpar{\fall{\Xila(\cA)} k}=O(C^k)$ as \ktoo, for some
$C<\infty$ depending on $\cA$ and $\gl$;
with only a little more effort, the same can be shown also for
$\cA\subseteq[a,\infty]$. 
This implies 
$\E e^{t\Xila(\cA)}<\infty$ for every such $\cA$ and $t<\infty$.
In particular, the distribution of $\Xila(\cA)$ is determined by its
(factorial) moments. Hence the formula in \refC{Cmom} in principle
determines the distribution of $\Xila(\cA)$ for any relatively compact
$\cA\subset(0,\infty]$.
Moreover, the formula in \refC{Cmom} easily extends to mixed factorial
moments of $\Xila(\cA_1),\dots,\Xila(\cA_m)$ when $\cA_1,\dots,\cA_m$ are
disjoint relatively compact Borel sets. This extension, which we leave
to the reader, characterizes the joint distribution 
of $\Xila(\cA_1),\dots,\Xila(\cA_m)$, and thus \cite[Theorem 3.1]{Kallenberg0}
the distribution of $\Xila$.
\end{remark}

\begin{remark}
  \label{Rmom2}
If $\cA\subseteq[a,\infty]$ for some $a>0$, we have, using the estimate
  in \refR{Rmom}, the standard formula
\begin{equation}
\P\bigpar{\Xila(\cA)=0}
=\sum_{k=0}^\infty\frac{(-1)^k}{k!}   \E\bigpar{\fall{\Xila(\cA)} k},
\end{equation}
which together with \refC{Cmom} (and perhaps the Bonferroni
inequalities) can be used for numerical evaluation
of $\P\bigpar{\Xila(\cA)=0}$, and thus, in particular, of the density
function in \refC{Cxi1}. 
\end{remark}

\begin{remark}
  \label{Rpalm*}
It follows easily from \refT{Tcomp} that a result analogous to
\refT{Tpalm} holds for $\Xix$ too. 
\end{remark}

\section{Limits as $\gl\to\pm\infty$} \label{Slimits}

In this section we consider limit results for $\Xila$, and in particular for
the largest point $\xila_1$, as $\gl\to\pm\infty$.
These results are equivalent to limit results for \gnp{} with
$p=n\qi+\gl(n)n^{-4/3}$ with $\gl(n)\to\pm\infty$ slowly, but we get
in this way no information on the allowed range of $\gl(n)$.

Consider first $\gl\to-\infty$.
By \eqref{f}, 
$F(x,\gl)\to\infty$ for every $x>0$ and
$F(x,\gl)$ is monotone in $\gl$ for $\gl\le0$.
Recalling the notation $\we\=\Xila[\eps,\infty)$,
it follows by
dominated convergence that, for every fixed $\eps>0$, 
$\E\we=\inteo\gL\gll(x)\dd x\to0$.
Hence, $\P(\we>0)\to0$ and $\P(\Xi[\eps,\infty)=\emptyset)\to1$.
Consequently, $\Xi\pto\emptyset$ (in the vague topology, see \refA{S:point})
and $\xila_1\pto0$.

We can by much more precise.
For $|\gl|>1$, let
\begin{equation}\label{al}
  \al\=3\ln|\gl|-\tfrac52\ln\ln|\gl|-\tfrac12\ln(2^43^5\pi),
\end{equation}
so that, as $\gl\to\pm\infty$, $\al\sim3\ln|\gl|$ and
\begin{equation}\label{q0}
e^{-\al}=|\gl|^{-3}(\ln|\gl|)^{5/2}(2^43^5\pi)\qq
\sim 4\pi\qq |\gl|^{-3}\al^{5/2}.
\end{equation}

\begin{theorem}\label{Tlim-}
  As $\gl\to-\infty$, 
  \begin{equation*}
\frac{|\gl|^2}2\xila_1-\al\dto V,	
  \end{equation*}
where $V$ has the Gumbel (extreme value) distribution
$\P(V\le s)=e^{-e^{-s}}$.
\end{theorem}
\begin{proof}
Fix a real $s$, and 
let $N\gll(x)\=\Xila(x,\infty)$, the number of points in $\Xila$
larger than $x$.
Thus $\E N\gll(x)=\E\Xila(x,\infty)=\int_x^\infty\gL\gll(y)\dd y$.
With the change of variables $y=2\gl\qii(\al+t)$, we obtain
\begin{equation}
  \label{q1}
\E N\gll\bigpar{2\gl\qii(\al+s)}
=\int_s^\infty2\gl\qii\gL\gll\bigpar{2\gl\qii(\al+t)}\dd t.
\end{equation}
For $\gl\le0$ and any real $t$ we have, by \eqref{f}, 
\begin{equation*}
F\bigpar{2\gl\qii(\al+t),\gl}
=\tfrac86\gl^{-6}(\al+t)^3 + 2|\gl|^{-3}(\al+t)^2+(\al+t)
=\al+t+o(1),
\end{equation*}
as $\gl\to-\infty$ with $t$ fixed.
Since $\Psi(x)\to1$ as $x\to0$, 
it follows from this,
\eqref{tint}
and \eqref{q0} that
\begin{multline*}
  2\gl\qii\gL\gll\bigpar{2\gl\qii(\al+t)}
\\
\begin{aligned}  
&
=(2\pi)\qqi2^{-3/2}|\gl|^3(\al+t)^{-5/2}
  \Psi\bigpar{2^{3/2}|\gl|^{-3}(\al+t)^{3/2}}
  e^{-F(2\gl\qii(\al+t),\gl)}
\\&
=\tfrac14\pi\qqi|\gl|^3\al^{-5/2}
e^{-\al-t+o(1)}
\to e^{-t}.
\end{aligned}
\end{multline*}
Moreover, for $t\ge s$ and $\gl<0$ with $|\gl|$ so large that
$\al>2|s|$ we also obtain, 
using $\Psi(x)=O\bigpar{e^{x^2/6}}$ from \eqref{jesper}
and $F(x,\gl)\ge x^3/6+x\gl^2/2$ from \eqref{f},
\begin{equation*}
  2\gl\qii\gL\gll\bigpar{2\gl\qii(\al+t)}
=O\Bigpar{|\gl|^3\al^{-5/2}
  e^{-(\al+t)}}
=O\bigpar{ e^{-t}}.
\end{equation*}
Consequently, we can use dominated convergence in \eqref{q1} and thus
\begin{equation*}
\E N\gll\bigpar{2\gl\qii(\al+s)}
\to
\int_s^\infty e^{-t}\dd t =e^{-s}.
\end{equation*}

Higher factorial moments can be computed similarly using \refC{Cmom},
with $B=(2\gl\qii(\al+s),\infty)$ and $x_j=2\gl\qii(\al+t_j)$.
Note that, for fixed $t_j$, $x_j\to0$, and thus 
$F(x_j,\gl-x_1-\dots-x_{j-1})=F(x_j,\gl)+o(1)$. 
Note further that, for $\gl<0$ and every $u,x\ge0$, 
$\gL\gxx{\gl-u}(x)\le\gL\gll(x)$; hence the bound used to verify
dominated convergence above applies to each factor in this
multivariate setting too. Consequently, for every $k\ge1$,
\begin{equation*}
\E \fall{\bigpar{N\gll\xpar{2\gl\qii(\al+s)}}}k
\to
\int_s^\infty \dotsm\int_s^\infty
e^{-t_1}\dotsm e^{-t_k} \dd t_k\dotsm\dd t_1 =\bigpar{e^{-s}}^k.
\end{equation*}
By the method of moments, this implies
$N\gll\bigpar{2\gl\qii(\al+s)}\dto\Po\bigpar{e^{-s}}$, and thus
\begin{equation*}
  \P\bigpar{\xila_1\le 2\gl\qii(\al+s)}
=\P\bigpar{N\gll\bigpar{2\gl\qii(\al+s)}=0}
\to e^{-e^{-s}}.
\end{equation*}
\vskip-\baselineskip
\end{proof}

\begin{remark}\label{Rlim-}
The proof yields also the asymptotic distribution of 
$\xila_2$, $\xila_3$, \dots. 
In fact, for every fixed $i$, as $\gl\to-\infty$,
\begin{equation*}
  \P\bigpar{\xilai\le 2\gl\qii(\al+s)}
=\P\bigpar{N\gll\bigpar{2\gl\qii(\al+s)}<i}
\to \sum_{j=0}^{i-1} \frac{e^{-js}}{j!}e^{-e^{-s}};
\end{equation*}
if we write the \rhs{} as $\P(V_i\le s)$, this can be written
\begin{equation*}
\frac{|\gl|^2}2\xila_i-\al\dto V_i.	  
\end{equation*}
Note that these asymptotic distributions are the same as for the
$i$:th records of suitable \iid{} sequences, see \cite[Section 2.2]{LLR}.
  
More generally, the proof above is easily extended to show that
$\Xila$ with the points rescaled as above, converges in distribution
to a Poisson process on $(-\infty,\infty)$ with intensity $e^{-s}$.
(This holds in $\fN[-a,\infty]$ for every $a$, say; we cannot use
$\fN(-\infty,\infty]$ directly, since the rescaled processes are not
elements of this space.)
Thus for $\gl\to-\infty$, the point process $\Xila$ becomes Poisson-like.
\end{remark}

In particular, $\xila_i$ is roughly $2\gl\qii\al\sim 6\ln|\gl|/\gl^2$
for every fixed $i\ge1$. This can be made precise in the following
form, where we use the notation that $X_\gl\psim x_\gl$ if $X_\gl/x_\gl\pto1$.

\begin{corollary}\label{Clim-}
As $\gl\to-\infty$,
$\xila_i\psim 6\ln|\gl|/\gl^2$ 
for every fixed $i\ge1$.
\nopf
\end{corollary}

\begin{remark}
The asymptotic results for large negative $\gl$ in
\refT{Tlim-} and \refR{Rlim-}
have a natural interpretation.  
The results on
the asymptotic distribution of $\xila_i$ are what they would be if $\Xila$ were
replaced by a Poisson point process with intensity $\gL\gll$.  This
corresponds to the view that as one moves in the critical window
toward the subcritical phase the largest components become ``local
phenomenon" and their interaction becomes negligible.
\end{remark}

Let us now turn to $\gl\to+\infty$. It is well-known that in this
case, with probability tending to 1, 
$\Xila$ contains exactly one large point.
In fact, $\xila_1\pto\infty$ and $\xila_2\pto0$ as \ltoox+.
Again, we can be much more precise.

Let $X$ and $Y$ by two random variables.
The \emph{total variation distance} between the distributions of $X$
and $Y$ is
defined as 
\begin{equation*}
  \dtv(X,Y)\=
\sup_B|\P(X\in B)-\P(Y\in B)|,
\end{equation*}
taking the supremum over all Borel sets $B$. 
Note that this only depends on the distributions $\cL(X)$ and $\cL(Y)$,
although we for simplicity use the notation $\dtv(X,Y)$ instead of 
$\dtv(\cL(X),\cL(Y))$; 
we will also write
$\dtv(X,\mu)$ when $Y$ has distribution $\mu$.
Note also that $\dtv$ is a very strong measure of distance between
distributions; for example, for a sequence $X_n$,
$\dtv(X_n,Y)\to0$ is much stronger than
$X_n\dto Y$, and thus (i) below is stronger than asymptotic normality in
the standard form \newline
$(\xila_1-2\gl)/\sqrt{2/\gl}\dto N(0,1)$.

\begin{theorem}\label{Tlim+}
If $\gl\to+\infty$, then
  \begin{romenumerate}
\item 
$\dtv\bigpar{\xila_1,N(2\gl,2\gl\qi)}\to0$;
\item
 $
\frac{\gl^2}2\xila_2-\al\dto V,	
 $
with $\al$ as in \eqref{al} and\/ $V$ as in \refT{Tlim-}.
  \end{romenumerate}
\end{theorem}

The proof below also shows that 
 $\tfrac12{|\gl|^2}\xila_i-\al\dto V_{i-1}$ for every $i\ge2$, with
 $V_i$ as in \refR{Rlim-}.	

\begin{corollary}\label{Clim+}
As $\gl\to+\infty$,
$\xila_1\psim 2\gl$ and
$\xila_i\psim 6\ln\gl/\gl^2$ 
for every fixed $i\ge2$.
\nopf
\end{corollary}

To prove \refT{Tlim+}, we begin with two lemmas.
Let $\phil$ denote the density
function of $N(2\gl,2\gl\qi)$; thus,
$\phil(x)=(\gl/4\pi)\qq e^{-\gl(x-2\gl)^2/4}$.

\begin{lemma}
  \label{LQ1}
As \ltoox+,
\begin{equation*}
\int_\gl^\infty \bigl|\gL\gll(x)-\phil(x) \bigr|\dd x\to0.  
\end{equation*}
\end{lemma}

The lower limit $\gl$ is for convenience only; it can easily be
replaced by, \eg,~1.

\begin{proof}
  For $x\ge\gl$, \eqref{jesper} yields
$\Psi(x^{3/2})=\tfrac12 x^3e^{x^3/24}\bigpar{1+o(1)}$, with $o(1)\to0$
as \ltoox+, uniformly in $x\ge\gl$.
Using $|a\qq-b\qq|=|a-b|/(a\qq+b\qq)$ and $|e^a-e^b|\le |a-b|e^{\max\set{a,b}}$,
we thus find from \eqref{tint} and \eqref{magnus}, for $x\ge\gl$,
\begin{align*}
  \gL\gll(x)
&=
(2\pi)\qqi \tfrac12 x\qq e^{-x(x-2\gl)^2/8}\bigpar{1+o(1)}
\\&
=
(8\pi)\qqi (2\gl)\qq e^{-x(x-2\gl)^2/8}\bigpar{1+o(1)}
+
O\Bigpar{|x-2\gl|\gl\qqi e^{-\gl(x-2\gl)^2/8} }
\\&
=\phil(x)\bigpar{1+o(1)}+
O\Bigpar{\bigpar{\gl\qq|x-2\gl|^3+\gl\qqi|x-2\gl|} e^{-\gl(x-2\gl)^2/8} }.
\end{align*}
The result follows by integrating; the $O$ term yields, 
if we let $Z\sim N(2\gl,4\gl\qi)$, 
$O\bigpar{\E (\gl\qq|Z-2\gl|^3+\gl\qqi|Z-2\gl|)}=O(\gl\qi)$.
\end{proof}

\begin{lemma}
  \label{Ldtv}
For any random variables $X$ and $Y$ with density functions $f_X$ and
$f_Y$, and any Borel set $B\subseteq\bbR$,
\begin{equation*}
  \dtv(X,Y)\le \int_B|f_X(x)-f_Y(x)|\dd x + \P(Y\notin B).
\end{equation*}
\end{lemma}

\begin{proof}
  It is well-known, and easy to verify, that 
$$
\dtv(X,Y)=\tfrac12\int_{-\infty}^\infty|f_X(x)-f_Y(x)|\dd x.
$$
Since $\int f_X=1=\int f_Y$, we have 
{\multlinegap=3pt
\begin{multline*}
  \int_{B^c}|f_X(x)-f_Y(x)|\dd x
\le
\int_{B^c}\bigpar{f_X(x)+f_Y(x)}\dd x
\\
=
2\P(Y\notin B) + \int_{B^c}\bigpar{f_X(x)-f_Y(x)}\dd x
=
2\P(Y\notin B) - \int_{B}\bigpar{f_X(x)-f_Y(x)}\dd x
\end{multline*}}
and thus
$
 \int_{\bbR}|f_X(x)-f_Y(x)|\dd x
\le
2\P(Y\notin B) + 2\int_{B}|f_X(x)-f_Y(x)|\dd x.
$
\end{proof}

\begin{proof}[Proof of \refT{Tlim+}]
 
If $x\ge\gl\ge0$, then $\gL\gllx(y)\le\gL\gxx0(y)$ for $y\ge0$ and thus
\begin{align*}
  \P\bigpar{\Xi\gllx(\gl,\infty)\ge1}
&
\le
  \E\Xi\gllx(\gl,\infty)
=
\int_\gl^\infty\gL\gllx(y)\dd y
\\&
\le
\int_\gl^\infty\gL\gxx0(y)\dd y
=
  \E\Xi\gxx0(\gl,\infty)
\to0
\end{align*}
as $\gl\to\infty$.
Hence, by \refC{Cxi1}, $\xila_1$ has a density function $\hla1$ with
$\hla1(x)=\bigpar{1-o(1)}\gL\gll(x)$ as $\gl\to\infty$, uniformly in $x\ge\gl$.
Since \refL{LQ1} implies 
$\int_\gl^\infty\gL\gll=O(1)$, this yields
\begin{equation*}
  \int_\gl^\infty \bigl|\hla1(x)-\gL\gll(x) \bigr|\dd x
\le   \E\Xi\gxx0(\gl,\infty) \int_\gl^\infty \gL\gll(x)\dd x
\to0.
\end{equation*}
Hence \refL{LQ1} yields 
$\int_\gl^\infty \bigl|\hla1(x)-\phil(x) \bigr|\dd x\to0$
as $\gl\to\infty$, and (i) follows by \refL{Ldtv}, with
$B=(\gl,\infty)$.

\newcommand\condd[2]{\cL\bigpar{#1\mid#2}}

For (ii), we observe that, by a simple extension of the proof of
\refC{Cxi1},
the conditional distribution $\cL(\Xila\mid\xila_1=x)$ 
equals the conditional distribution 
$\condd{\gd_x+\Xi\gllx}{\Xi\gllx(x,\infty)=0}$.
Since the second largest point in 
$\gd_x+\Xi\gllx$, when $\Xi\gllx(x,\infty)=0$, is the largest point
$\xi\gllx_1$ in $\Xi\gllx$, we have, in particular,
\begin{equation}
  \label{q3}
\condd{\xila_2}{\xila_1=x}
=
\condd{\xi\gllx_1}{\xi\gllx_1\le x}.
\end{equation}

Let \ltoox+, and assume $\gl>2$.
By (i), $\P(|\xila_1-2\gl|<1)\to1$. 
If $|x-2\gl|<1$ and $\gl'\=\gl-x$, then $|\gl'-(-\gl)|<1$ and thus, by
\eqref{al}, $|\al-\alx|=O(1/\gl)$. Hence, it follows from \refT{Tlim-}
that 
\begin{equation}\label{q6}
  \tfrac12{\gl^2}\xi\gllx_1-\al
=|\gl/\gl'|^2 \bigpar{\tfrac12{|\gl'|^2}\xi\gxx{\gl'}_1-\alx}
+\bigpar{|\gl/\gl'|^2-1}\alx
+\alx-\al
\dto V.
\end{equation}
Furthermore, 
$\P\bigpar{\xi\gllx_1\le x}\to1$, again by \refT{Tlim-}, and thus
\eqref{q6} holds also for the conditional distribution given
$\xi\gllx_1\le x$.
By \eqref{q3} and $\P(|\xila_1-2\gl|<1)\to1$,
this yields, for every $y$,
\begin{equation*}
  \P\bigpar{\tfrac12\gl^2\xila_2-\al\le y}
=\E  \P\bigpar{\tfrac12\gl^2\xila_2-\al\le y\mid\xila_1}
\to\P(V\le y),
\end{equation*}
which proves (ii).
\end{proof}

\begin{remark}
  Note that it is a fallacy to believe that \refT{Tpalm} implies that
  $\Xila$ conditioned on $\xila_1=x$ has the distribution of
  $\gd_x+\Xi\gllx$; as is seen in the proof above, the correct
  conclusion requires conditioning on $\Xi\gllx(x,\infty)=0$.
Nevertheless, the proof also shows that the erroneous statement is
  asymptotically correct as \ltoox+:
If $\xi$ has the distribution of $\xila_1$ given in \refC{Cxi1}, or
  simply $\xi\sim N(2\gl,2\gl\qi)$, and given $\xi$ we take a random
  $\Xi\gxx{\gl-\xi}$, then the distribution of 
$\gd_\xi+\Xi\gxx{\gl-\xi}$ approximates
that of
$\Xila$, and 
$\dtv\bigpar{\gd_\xi+\Xi\gxx{\gl-\xi},\Xila}\to0$ as \ltoox+.
\end{remark}

\begin{remark}
As remarked above, \refT{Tlim+} implies asymptotic normality of the 
size of the largest  component in \gnp{} 
in the case $p=n\qi+\gl(n)n^{-4/3}$ with $\gl(n)\to\infty$ slowly
(without specifying the allowed rate).
Indeed, asymptotic normality has been shown for all $\gl(n)$ 
in the range $\gl(n)\to\infty$  but $\gl(n)=O(n\qw)$, \ie{} $p=O(n\qi)$,
by   Pittel \cite{Pittel} ($p=c/n$) and
Pittel and Wormald \cite{PittelW} (the general case).
\end{remark}

\appendix
\section{Appendix: Point processes}\label{S:point}
We give here some technical remarks on point processes;
see \eg{} \cite{Kallenberg0} and
\cite[Section 4]{SJ136} for further details and proofs.

Let $\fS$ be a `nice' topological space 
(more precisely, a locally compact Polish space);
in this paper we only consider the intervals 
$(0,\infty)$ and 
$(0,\infty]$ and 
their products with $\bbN$ or $\bbNx$.
Although we regard a point process as a random (multi)set
$\set{\xi_i}_i\subset\fS$, it is technically convenient to formally
define it as a random measure
$\sum_i\delta_{\xi_i}$.
Hence, if $\Xi$ denotes the point process \set{\xi_i}, we write $\Xi(A)$ for 
the number of points $\xi_i$ that belong to a subset $A\subseteq\fS$;
similarly,
for suitable functions  $f$ on $\fS$,
$\int f\,d\Xi=\sum_i f(\xi_i)$.

Thus, let $\fN=\fN(\fS)$ be the class of all Borel measures $\mu$ on $\fS$ such
that
$\mu(A)$ is a (finite) integer $0,1,\dots$ for every relatively
compact Borel set $A$; this coincides with the class of all finite or
countably infinite sums
of the type 
$\sum_i\delta_{x_i}$, where $x_i\in\fS$ and 
each  compact subset of $\fS$ contains
only a finite
number of $x_i$,
and we identify such a sum with the (multi)set $\set{x_i}$.

The standard topology on $\fN$ (known as the \emph{vague topology})
is defined such that,
for $\mu,\mu_1,\mu_2,\dots\in\fN$, 
$\mu_n\to\mu$ if and only if $\int f\,d\mu_n\to\int f\,d\mu$ for
every $f\in\Cc(\fS)$, the space of
(real-valued) continuous functions on $\fS$ with compact support.
(This is a metrizable topology and $\fN$ is a Polish space, see
\cite[Section 15.7]{Kallenberg0}.)

A point process on $\fS$ is a random element of $\fN$.
If $\Xi$ is a point process on $\fS$, there exists a
unique Borel measure $\nu$ on $\fS$ such that $\E\Xi(A)=\nu(A)$ for
every Borel set $A$, and more generally $\E\int h\,d\Xi=\int h\,d\nu$
for every positive measurable function $h$. This measure $\nu$ is called
the \emph{intensity} of $\Xi$. In the cases we consider, $\fS$ is an
interval or a union of intervals,
and $\nu$ is absolutely continuous; then also the function $d\nu/dx$
is called the intensity.

If $\Xi_n$ and $\Xi$ are point processes on $\fS$, then
$\Xi_n\dto\Xi$ (w.r.t.\ the vague topology just defined)
if and only if $\int f\,d\Xi_n\dto\int f\,d\Xi$ (as real-valued random
variables) for every $f\in\Cc(\fS)$.
It is also true that
$\Xi_n\dto\Xi$
if and only if $\Xi_n(A)\dto\Xi(A)$ for every relatively compact Borel
set $A\subseteq\fS$ such that $\Xi(\partial A)=0$ a.s., and moreover
joint convergence 
holds for every finite collection of such sets $A$.

We state a particular case that we need.
Say that a point $x$ is a \emph{continuity point} of a point process
$\Xi$ 
if $x$ is a continuity point of $\E\Xi$,
\ie{} if
$\E\Xi\set{x}=0$, or equivalently,
$x\notin\Xi$ a.s.

\begin{lemma}\label{L:point3}
  If\/ $\Xi_n\dto\Xi$ as point processes on an interval $J$, then
  $\Xi_n[a,b]\dto\Xi[a,b]$ for every interval $[a,b]\subset J$ such that
$a$ and $b$ are continuity points of $\Xi$.
\nopf
\end{lemma}

Note that the definitions of both point processes and convergence of
them are sensitive to the choice of $\fS$, 
since a point process is not allowed to have any cluster point in $\fS$.
Hence, it matters
whether 
boundary points are included in $\fS$, even if they are not attained
by any point.
For example, if $\fS$ is the closed
interval $[0,\infty]$ (or any compact set), then every point process is finite.
If, instead, $\fS$ is the half-open interval $(0,\infty]$, 
then an element $\mu\in\fN$
is finite on every interval $[a,\infty]$, and thus every point
process
may be written as a (finite or infinite) set $\set{\xi_i}$ with
$\infty\ge\xi_1\ge\xi_2\ge\dots$ and, if the set is infinite, $\xi_i\to 0$ as
$i\to\infty$.
Similarly, a point process on the open interval $(0,\infty)$ may have both $0$
and $\infty$ as cluster points.
By including one or both endpoints,
we thus get stronger conditions, and,
similarly, we get a stronger mode of convergence.
It may thus be advantageous to consider (when possible)
a set of points in $(0,\infty)$  as a point process on $[0,\infty)$, 
$(0,\infty]$ or $[0,\infty]$.

For point processes on a closed or half-open interval, with the points
ordered as above,  
convergence 
is equivalent to joint convergence of the
individual points. We state this for the case we are interested
in. 

\begin{lemma}\label{L:point2}
There is a bijection between $\fN(0,\infty]$ and the space of
sequences $(\xi_i)_1^\infty$ with $\xi_1\ge\xi_2\ge\dots\ge0$ and
$\lim_{i\to\infty} \xi_i=0$,
such that $\Xi=\set{\xi_{i}}_{i=1}^{N}\in\fN$ (or, more formally,
$\Xi=\sum_{i=1}^{N}\gd_{\xi_{i}}$),
with
$\xi_{1}\ge\xi_{2}\ge\dots$ and $0\le N\le\infty$, corresponds to the sequence
$(\xi_i)_1^\infty$ where we define $\xi_i\=0$ for $i>N$.
This bijection is a homeomorphism between $\fN$ with the vague
topology and the space of sequences with component-wise convergence
(\ie, the restriction of the product topology on $[0,\infty]^\infty$).

Consequently, if $\Xi_n$, $1\le n\le\infty$, are point processes on the
interval $(0,\infty]$,
and we write $\Xi_n=\set{\xi_{ni}}_{i=1}^{N_n}$ with
$\xi_{n1}\ge\xi_{n2}\ge\dots$ and $0\le N_n\le\infty$, and
if some $N_n<\infty$, we further define $\xi_{ni}=0$
for $i>N_n$,
then $\Xi_n\dto\Xi_\infty$ if and only if
$
(\xi_{n1},\xi_{n2},\dots)
\dto(\xi_{\infty1},\xi_{\infty2},\dots)$,
in the standard sense that all finite dimensional distributions
converge.
\nopf
\end{lemma}

\end{document}